\documentclass{ps}
\usepackage{layout, graphicx, amsmath, amsthm, amssymb, euscript, enumerate, bbold, bbm, graphicx, fancyhdr,xcolor,enumitem}
\numberwithin{equation}{section}
\usepackage{tikz}
\usepackage{yfonts} 
\usepackage[toc,page]{appendix}

\usepackage[pdftex,bookmarks,colorlinks,pdfstartview={FitH}]{hyperref}
\usepackage{subcaption}

\hypersetup{colorlinks=false}

\def\I{\mathfrak{I}}

\def\B{\mathcal{B}}

\def\N{\mathbb{N}}

\def\P{\mathbb{P}}

\def\R{\mathbb{R}}
\def\E{\mathbb{E}}

\def\eg{\textit{e.g.} }
\def\ie{\textit{i.e.} }

\def\Var{\mbox{Var}}

\theoremstyle{definition}
\newtheorem{assmptn}[thrm]{Assumption}

\begin{document}

	\title{Chain-referral sampling on Stochastic Block Models}\thanks{This work was done during the PhD thesis of the author under the supervision of Jean-St\'ephane Dhersin and Tran Viet Chi. The author was partially supported by the Chaire MMB (Mod\'elisation Math\'ematique et Biodiversit\'e of Veolia-Ecole Polytechnique-Museum National d'Histoire Naturelle-Fondation X) and by the ANR Econet (ANR-18-CE02-0010).}
	
\author{Thi Phuong Thuy Vo}\address{Vo Thi Phuong Thuy, Univ. Paris 13, CNRS, UMR 7539 - LAGA, 99 avenue J.-B. Cl\'ement, F-93430 Villetaneuse, France; E-mail: phuongthuywz@gmail.com}

\date{\today}
	\begin{abstract} 
		The discovery of the ``hidden population", whose size and membership are unknown, is made possible by assuming that its members are connected in a social network by their relationships. We explore these groups by a chain-referral sampling (CRS) method, where participants recommend the people they know. This leads to the study of a Markov chain on a random graph where vertices represent individuals and edges connecting any two nodes describe the relationships between corresponding people. We  are interested in the study of CRS process on the stochastic block model (SBM), which extends the well-known Erd\"os-R\'enyi graphs to populations partitioned into communities. The SBM considered here is characterized by a number of vertices $N$, a number of communities (blocks) $m$, proportion of each community $\pi=(\pi_1,...,\pi_m)$ and a pattern for connection between blocks $P=(\lambda_{kl}/N)_{(k,l) \in \{1,...,m\}^2}$. In this paper, we give a precise description of the dynamic of CRS process in discrete time on an SBM. The difficulty lies in handling the heterogeneity of the graph. We prove that when the population's size is large, the normalized stochastic process of the referral chain behaves like a deterministic curve which is the unique solution of a system of ODEs.
	\end{abstract}

\begin{resume}Nous nous int\'eressons à l'\'etude de ``populations cach\'ees'', de taille inconnue, et dont on ne connaît pas les membres. La d\'ecouverte d'une population cach\'ee est rendue possible en supposant que ses individus sont connect\'es par un r\'eseau social. Nous explorons ces groupes par une m\'ethode de sondage par cha\^inage (``Chain referral sampling", CRS), o\`u les r\'epondants recommandent leurs contacts. Ceci conduit à l'\'etude d'une cha\^ine de Markov sur un graphe al\'eatoire dont les sommets repr\'esentent les individus et dont les ar\^etes d\'ecrivent les relations entre les deux personnes qu'elles relient. Les personnes interrog\'ees sont invit\'ees \`a indiquer leurs partenaires et un certain nombre de coupons est remis \`a certaines de ces personnes. Le sondage par cha\^inage recherche les noeuds cach\'es dans la population en suivant au hasard les ar\^etes du r\'eseau social sous-jacent, ce qui permet de tracer les individus \'echantillonn\'es.  Nous \'etudions le processus CRS lorsque le r\'eseau est un mod\`ele à blocs stochastiques (``Stochastic Block Model", SBM), qui est une extension du mod\`ele d'Erd\"os-R\'enyi aux populations partitionn\'ees en communaut\'es. Le SBM consid\'er\'e ici est caract\'eris\'e par un certain nombre de sommets $ N $ (taille de la population), un certain nombre de communaut\'es (blocs) $ m $, une distribution de blocs $ \pi = (\pi_1, ..., \pi_m) $ repr\'esentant la proportion de chaque communaut\'e et une matrice permettant de d\'efinir les liens entre sommets appartenant \`a des blocs donn\'es  $ P = (\lambda_ {kl} / N)_{(k, l) \in \{1, ..., m \} ^ 2} $. Dans cet article, nous donnons une description pr\'ecise de la dynamique du processus CRS en temps discret sur un SBM. La difficult\'e r\'eside dans la gestion de l'h\'et\'erog\'en\'eit\'e du graphe. Dans notre mod\`ele, le graphe et la marche al\'eatoire sont construits simultan\'ement. Ensuite, nous \'etudions l'\'evolution de cette cha\^ine en consid\'erant le processus normalis\'e sur l'\'echelle de temps $ [0,1] $. Nous d\'{e}montrons que lorsque la taille de la population est grande, le processus aléatoire CRS normalis\'e se comporte comme une courbe d\'eterministe qui est la solution unique d'un syst\`eme d'ODE.
\end{resume}


\maketitle
	
	\section{Introduction}

	In Sociology, some populations may be hidden because their members share common attributes that are illegal or stigmatized. These hidden groups may be hard to approach because these individuals try to conceal their identities due to legal authorities (e.g. drugs users) or because of the social pressure (e.g. men having sex with men). In such populations, all the information is unknown: there is no sampling frame such as lists of the members of the population or of the relationship between the latter. It causes many challenges for researchers to identify these groups. The discovery of the hidden populations is made possible by assuming that its members are connected by a social network. The population is described by a graph (network) where each individual is represented by a vertex and any interaction or relationship (e.g. friendship, partnership) between a couple of individuals is represented by an edge matching the corresponding vertices. Thanks to this important feature, we are allowed to investigate these populations by using a Chain-referral Sampling (CRS) technique, such as snowball sampling, targeting sampling, respondent driven sampling etc. (see the review of \cite{Shaghaghi} or \cite{goodman1961,heckathorn,heckathorn2}). CRS consists in detecting hidden individuals in a population structured as a random graph, which is modeled by a stochastic process that we study here.  The principle of CRS is that from a group of initially recruited individuals, we follow their connections in the social network to recruit the subsequent participants. The exploration proceeds from node to node along the edges of the graph.  The interviewees induce a sub-tree of the underlying real graph, and the information coming from the interviews gives knowledge on other non-interviewed individuals and edges, providing a larger sub-graph. We aim at understanding this recruitment process from the properties of the explored random graph. The CRS showed its practicality and efficiency in recruiting a diverse sample of drug users (see \cite{saadati}).\\

	CRS models are hard to study from a theoretical point of view without any assumption on the graph structure. In this paper, we consider a particular model with latent community structure: the stochastic block model (SBM) proposed by {Holland et al.}\cite{Holland}. This model is a useful benchmark for some statistical tasks as recovering community (also called blocks or types in the sequel) structure in network science \cite{Gadde,Girvan,barbillon}. By block structure, we mean that the set of vertices in the graph is partitioned into subsets called blocks and nodes connect to each other with probabilities that depend only on their types, \ie the blocks to which they belong. For example, edges may be more common {within a block} than between blocks (e.g. group of people having sexual contacts). We {recall here} the definition of SBM (we refer the reader to the survey in \cite{abbe}):
	\begin{dfntn}
		Let $N$ be a positive integer (number of vertices), $m$ be a positive integer (number of blocks or types), $\pi=(\pi_1,...,\pi_m)$ be a probability distribution on $\{1,\dots m\}$ (the probabilities of the $m$ types, \ie a vector of $[0,1]^m$ such that $\sum_{k=1}^m \pi_k=1$) and $P=(p_{kl})_{(k,l)\in \{1,...,m\}^2}$ be a symmetric matrix with entries $p_{kl}\in [0,1]$ (connectivity probabilities). The pair $(\Gamma, G)$ is drawn under the distribution SBM$(N,\pi,P)$ if the vector of types $\Gamma$ {is an $N$-dimensional random vector, whose components are i.i.d., $\{1, \ldots, m\}$-valued with the law $\pi$, and $G$ is} a simple graph of size $N$ where vertices $i$ and $j$ are connected independently of other pairs of vertices with probability $p_{\Gamma_i\Gamma_j}$. We also denote the blocks (community sets) by: $[l]: = \{v \in \{1,...,N\}: \Gamma_v = l\}$ with the size $N_l:=|[l]|, l \in \{1,...,m\}$.		
	\end{dfntn}
	Notice that when $m=1$, \ie there is only one type. Any arbitrary pair of vertices is connected independently to the others with the same probability $p_{11}$, SBM becomes the Erd\"os-R\'enyi graph, which is studied in \cite{vo1}. \\
	Here, we consider the Poisson case where the connectivity probabilities $p_{kl}$ depend on $N$ and are given by $p_{kl}=\lambda_{kl}/N$. This means that each individual of the block $k$ {contacts in average $\lambda_{kl} \pi_l$ individuals} of the block $l$. This implies that the network examined is sparse. In the present work, we give a rigorous description of a CRS on such SBM and study the propagation of {the referral chain} on this sparse model.\\

\begin{figure}[h!]
		\centering
		\begin{tabular}{|c c |c|}
			\hline
			\begin{tikzpicture}[line cap=round,line join=round,x=0.7cm,y=0.7cm]
			\draw[blue,fill] (0,0) circle (2.5pt);
			\draw[] (0.5,2)circle (2.5pt); \draw[] (2.5,0)circle (2.5pt); \draw[] (1.5,-1.2)circle (2.5pt); \draw[] (4,2.2)circle (2.5pt); \draw[] (5,1)circle (2.5pt); \draw[] (5,-1)circle (2.5pt); \draw[] (4,-2.2)circle (2.5pt); \draw[] (7,0.5)circle (2.5pt);
			\draw[] (6,2.3)circle (2.5pt);
			\draw[] (6.5,-1.2)circle (2.5pt);
			\draw[] (1.5,1.5)circle (2.5pt);
			\draw[] (0.5,-2)circle (2.5pt);
			\end{tikzpicture}
			&
			&
			\begin{tikzpicture}[line cap=round,line join=round,x=0.7cm,y=0.7cm]
			\draw[] (4,2.2)circle (2.5pt); \draw[] (5,1)circle (2.5pt); \draw[] (5,-1)circle (2.5pt); \draw[] (4,-2.2)circle (2.5pt); \draw[] (7,0.5)circle (2.5pt);
			\draw[] (6,2.3)circle (2.5pt);
			\draw[] (6.5,-1.2)circle (2.5pt);
			\draw[] (1.5,1.5)circle (2.5pt);
			\draw[] (0.5,-2)circle (2.5pt);
			
			\draw [red,fill] (0,0) circle (2.5pt);
			\draw[blue,fill] (0.5,2)circle (2.5pt); \draw[blue,fill] (2.5,0)circle (2.5pt); \draw[gray,fill] (1.5,-1.2)circle (2.5pt);
			\draw[line width=1.2pt] (0,0) -- (0.5,2); \draw[line width=1.2pt] (0,0) -- (2.5,0); \draw[line width=1.2pt] (0,0) -- (1.5,-1.2);
			\end{tikzpicture}
			\\
			& & \\
			Step 0 & & Step 1\\
			\hline
			& & \\
			\begin{tikzpicture}[line cap=round,line join=round,x=0.7cm,y=0.7cm]
			\draw[] (4,2.2)circle (2.5pt); \draw[] (5,1)circle (2.5pt); \draw[] (5,-1)circle (2.5pt); \draw[] (4,-2.2)circle (2.5pt); \draw[] (7,0.5)circle (2.5pt);
			\draw[] (6,2.3)circle (2.5pt);
			\draw[] (6.5,-1.2)circle (2.5pt);
			\draw[] (1.5,1.5)circle (2.5pt);
			\draw[] (0.5,-2)circle (2.5pt);
			
			\draw [red,fill] (0,0) circle (2.5pt);
			\draw[blue,fill] (0.5,2)circle (2.5pt);  \draw[gray,fill] (1.5,-1.2)circle (2.5pt);
			\draw[line width=1.2pt] (0,0) -- (0.5,2); \draw[line width=1.2pt] (0,0) -- (2.5,0); \draw[line width=1.2pt] (0,0) -- (1.5,-1.2);
			
			\draw[red,fill] (2.5,0) circle (2.5pt);  \draw[gray,fill] (4,2.2)circle (2.5pt); \draw[blue,fill] (5,1)circle (2.5pt); \draw[gray,fill] (5,-1)circle (2.5pt); \draw[blue,fill] (4,-2.2)circle (2.5pt);
			\draw[line width=1.2pt] (2.5,0) -- (4,2.2);
			\draw[line width=1.2pt] (2.5,0) -- (5,1);
			\draw[line width=1.2pt] (2.5,0) -- (5,-1);
			\draw[line width=1.2pt] (2.5,0) -- (4,-2.2);
			\end{tikzpicture}
			& &
			\begin{tikzpicture}
			[line cap=round,line join=round,x=0.7cm,y=0.7cm]
			\draw[] (4,2.2)circle (2.5pt); \draw[] (5,1)circle (2.5pt); \draw[] (4,-2.2)circle (2.5pt);
			\draw[] (6,2.3)circle (2.5pt);
			\draw[] (6.5,-1.2)circle (2.5pt);
			\draw[] (1.5,1.5)circle (2.5pt);
			\draw[] (0.5,-2)circle (2.5pt);
			
			\draw [red,fill] (0,0) circle (2.5pt);
			\draw[blue,fill] (0.5,2)circle (2.5pt);  \draw[gray,fill] (1.5,-1.2)circle (2.5pt);
			\draw[line width=1.2pt] (0,0) -- (0.5,2); \draw[line width=1.2pt] (0,0) -- (2.5,0); \draw[line width=1.2pt] (0,0) -- (1.5,-1.2);
			
			\draw[red,fill] (2.5,0) circle (2.5pt);  \draw[gray,fill] (4,2.2)circle (2.5pt); \draw[blue,fill] (4,-2.2)circle (2.5pt);
			\draw[line width=1.2pt] (2.5,0) -- (4,2.2);
			\draw[line width=1.2pt] (2.5,0) -- (5,1);
			\draw[line width=1.2pt] (2.5,0) -- (5,-1);
			\draw[line width=1.2pt] (2.5,0) -- (4,-2.2);
			
			\draw[red,fill] (5,1)circle (2.5pt);
			\draw[blue,fill] (7,0.5)circle (2.5pt);
			\draw[blue,fill] (5,-1)circle (2.5pt);
			\draw[line width=1.2pt] (5,1) -- (7,0.5);
			\draw[line width=1.2pt] (5,1) -- (5,-1);
			\end{tikzpicture}\\
			& & \\
			Step 2 & & Step 3\\
			\hline
		\end{tabular}
		\begin{tikzpicture}
		\draw[red,fill] (0.5,-0.5) circle (2.5pt) node [right] {\quad off-mode node (who has been interviewed)};
		\draw[blue,fill] (0.5,-1) circle (2.5pt) node [right] {\quad active node (who has coupon but has not been interviewed yet)};
		\draw[gray,fill] (0.5,-1.5) circle (2.5pt) node [right] {\quad explored but still inactive node (who has been named but did not receive coupons)};
		\end{tikzpicture}
		\caption{Description of how the chain-referral sampling works. In our model, the random network and the CRS are constructed simultaneously. For example at step 3, an edge between two vertices who are already known at step 2 is revealed.} \label{fig:description}
	\end{figure}
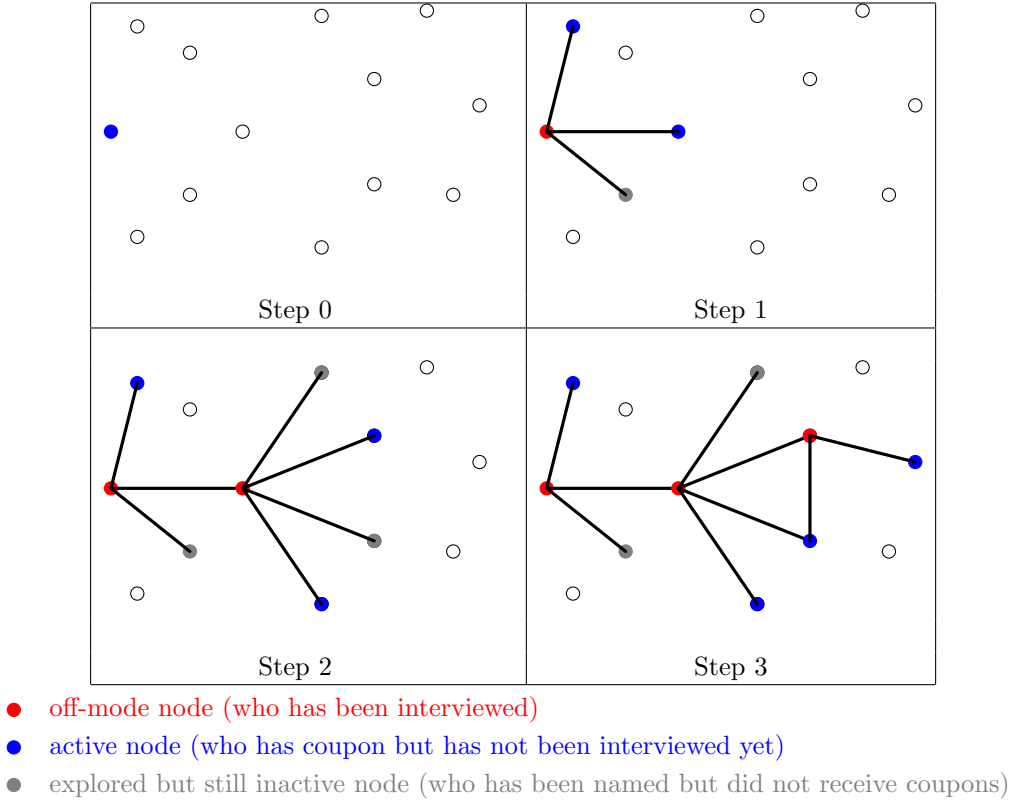

	The CRS relies on a random peer-recruitment process. To handle the two sources of randomness, the graph and the exploring process on it are constructed simultaneously. In the construction, the vertices of the graph will be in 3 different states: inactive vertices that have not being contacted for interviews, active vertices that constitute the next interviewees and off-mode vertices that have been already interviewed. {The idea to describe the random graph as a Markov exploration process with active, explored and unexplored nodes is classical in random graphs theory. It has been used as a convenient technique to expose the connections inside a cluster, especially to discover the giant component in a random graph models, for example see \cite{durrett2007, remcovanderhofstad}. In our case, there is a slight difference in the recruiting process: the number of nodes being switched to the active mode is set to be bounded by a constant. This trick helps to improve the bias towards high-degree nodes in the population (see \cite{heckathorn2}).} At the beginning of the survey, all individuals in the population are hidden and are marked as inactive vertices. We choose some people as seeds of the investigation and activate them. During the interview these individuals name their contacts and a maximum number $c$ of coupons are distributed to the latter, who become active nodes. One by one, every carrier of a coupon can come to a private interview and is asked in turn to give the names of her/his peers. Whenever a new person is named, one edge connecting the interviewee and her/his contact is added but they remain inactive until they receive a coupon. After finishing the interview, a maximum number of $c$ new contacts receive one coupon each and are activated. So if the interviewee names more than $c$ people, a number of them are not given any coupon and can be still explored later provided another interviewee mentions them. After that, the node associated to the person who has just been interviewed is switched to off-mode and is no longer recruited again, {see Figure} \ref{fig:description}. We repeat the procedure of interviewing, referring, distributing coupons until there is no more active vertex in the graph (no more coupon is returned). Each person returning a coupon receives some money as a reward for her/his participation, and an extra bonus depending on the number contacts that will later return the coupons. Notice that each individual in the population is interviewed just once and we assume here that there is no restriction on the total number of coupons.\\
	
	The process of interest counts the number of coupons present in the population. We also want to know how many people are detected, which leads to the number of people explored but without coupons. Denote by the discrete time $n \in \N=\{0,1,2,\dots \}$ the number of interviews completed, $A_n$ corresponds to the number of individuals that have received coupons but that have not been interviewed yet (number of active vertices); $B_n$ to the number of individuals cited in the interviews but who have not been given any coupon (number of found but still inactive vertices) and $U_n$ to the total number of individuals having been interviewed (number of off-mode nodes). \\
	Because of the connectivity properties of the SBM graphs, we need to keep track of the types of the interviewees and the coupons distributed not only to one community but also in general to each of the $m$ communities at every step. We then associate to the chain-referral the following stochastic vector process $X_n:=(A_n,B_n,U_n), \quad n \in \N$:
	\begin{equation*}
	X_n:= \begin{pmatrix}
	A_n \\ B_n \\U_n
	\end{pmatrix} =
	\begin{pmatrix}
	A_n^{(1)} & \cdots &A_n^{(m)}\\
	B_n^{(1)} & \cdots & B_n^{(m)}\\
	U_n^{(1)} & \cdots & U_n^{(m)}
	\end{pmatrix}, \quad n \in \mathbb{N},
	\end{equation*} where $A_n^{(l)}$ (resp. $B_n^{(l)}$ and $u_n^{(l)}$) corresponds to the number of active nodes (resp. of found but inactive nodes and of off-mode nodes) of type $l$ at step $n$. In all the paper, we will use the notation $(X_n^{1,(l)},X_n^{2,(l)},X_n^{3,(l)})=(A^{(l)}_n,B_n^{(l)},U_n^{(l)})$.

	The main object of the paper is to establish an approximation result when the size $N$ of the SBM graph tends to infinity. In this case, the chain-referral process correctly renormalized is:
	\begin{equation}\label{eq:renorm}
	X^N_t:=\frac{1}{N}X_{\lfloor Nt \rfloor}= \left(\frac{A_{\lfloor Nt \rfloor}}{N},\frac{B_{\lfloor Nt \rfloor}}{N},\frac{U_{\lfloor Nt \rfloor}}{N}\right)\in [0,1]^{3\times m},\quad t \in [0,1].
	\end{equation}

In all the paper, we consider spaces $\R^d$ equipped with the $L^1$-norm defined for $x=(x^1,...,x^d)$ as $\|x\|=\sum_{k=1}^{d}|x^k|$. For all $N$, the process $X_\cdot^N$ lives in the space of c\`adl\`ag processes $\mathcal{D}([0,1],[0,1]^{3\times m})$ {equipped with Skorokhod topology} (see \cite{ethierkurtz,jakubowski,joffemetivier}).

{There exist to our knowledge a few works of studying CRS form a probabilistic point of view, for example Athreya and Roellin \cite{athreya}. In their work, they obtained a result in a slightly different framework:}  they consider random walks on the limiting graphon to construct a sequence of sub-graphs, which converges almost surely to the graphon underlying the network in the cut-metric. Whereas we take here to the limit both the graph and its exploring random walk simultaneously. The main result of this paper is that the process $(X_.^N)_N$ converges to a system of ordinary differential equations (ODEs). {There has also been literature on random walks exploring graphs possibly with different mechanism (see \cite{bollobasriordan2012, enriquezfaraudmenard} for instance). Here we allow the exploring Markov process to branch. Also, our process bares similarities with epidemics spreading on graphs (see \cite{barbourreinert,brittonpardoux2019,decreusefond2010,jansonluczakwind2014}) but with the additional constraint of a maximum number of distributed coupons here.}\\

{The CRS is constructed by the similar principle of an epidemic spread and starts with a single individual. There are two main phases of evolution (see \cite{barbourreinert}): the initial phase is well approximated by a branching process (which we are neglecting here) and the second phase is when the stochastic process is approximated by an deterministic curve. In this paper, we focus on the second phase, but let us comment quickly on the first phase. In the sequel, we will assume that:}
	
	\begin{assmptn} \label{assp:transitionofsbm}
	{For each $\ell, k \in \{1, \ldots, m \}$, denote $\mu_{\ell k} = \lambda_{\ell k} \pi_k$. We assume that the matrix $\mu=(\mu_{\ell k})_{\ell, k \in \{1,..., m\}}$ is \textit{irreducible} and the largest eigenvalue of $\mu$ is larger than $1$.}
	\end{assmptn}

\begin{rmrk}
		Under the Assumption \ref{assp:transitionofsbm}, from the proof of Theorem 3.2 of Barbour and Reinert \cite{barbourreinert}, the early stages of the CRS is now can be associated approximated by a multitype branching process with the offspring distributions determined by the matrix $\mu$. {Thanks to the Assumption \ref{assp:transitionofsbm} the multitype branching process associated with the offspring matrix $\mu$ is supercritical.} The analogous results for the extinction probability and for the number of offspring at the $n^{th}$ generation as in the single branching process have been proved in Chapter 5 of \cite{athreyajagers}: the mean matrix of the population size at time $n$ is proportional to $\mu^n$. And follow the claim (3.11) of Barbour and Reinert \cite{barbourreinert}, we can deduce that if we start with a single individual, then after a finite steps, we can reach a positive fraction of explored individuals in the population with a positive probability.
\end{rmrk}

\begin{assmptn} \label{assp:initialcondition}
Set $a_0, b_0,u_0 \in [0,1]^m$, $a_0 = (a_0^{(1)}, \ldots, a_0^{(m)})$ such that $\sum_{i=1}^m a_0^{(i)} = \|a_0\| \in [0,1]$, and set $b_0, u_0\in [0,1]^m$, with $b_0= (0, \ldots, 0)$ and $u_0 =(0, \ldots, 0)$. We assume that the sequence $X_0^N = \frac{1}{N} X_0$ converges in probability to the vector $(a_0, b_0, u_0)$, as $N \rightarrow +\infty$.
\end{assmptn}
It means that the initial number of individuals with type $i$ at the beginning of the survey is approximately $\lfloor a_0^{(i)}N \rfloor$.	A possible way to initializing the process is to draw $A_0$ from a multinomial distribution $\mathcal{M}(\lfloor \|a_0\| N \rfloor; \pi_1, \ldots, \pi_m)$.\\

	\begin{thrm}\label{th:main} Under the assumptions \ref{assp:transitionofsbm} and \ref{assp:initialcondition}, we have: when $N$ tends to infinity, the process $(X_\cdot^N)_N$ converges in distribution in $\mathcal{D}([0,1],[0,1]^{3\times m})$ to a deterministic vectorial function $x=(x_\cdot^{(l)})_{1\leq l \leq m}=(a^{(l)}_\cdot,b^{(l)}_\cdot,u^{(l)}_\cdot)_{1 \leq l \leq m}$ in $\mathcal{C}([0,1],[0,1]^{3\times m})$, which is the unique solution {of the system of differential equations}
		\begin{equation}\label{limiting_eq}
		x_t=x_0 + \int_{0}^{t} f(x_s)ds,
		\end{equation}where
		  $f(x_s):= (f_{il}(x_s))_{\substack{1\leq i \leq 3 \\1\leq l \leq m}}$ has an explicit formula described as follows. Denote
		 \begin{equation} t_0 := \inf\{ t \in [0,1]: \| a_t\| := a^{(1)}_t + \ldots + a^{(m)}_t=0\}. \end{equation}		
		 For $s \in [0,t_0]$,		
		\begin{align}
		f_{1l}(x_s) &= \sum_{k=1}^{m} \dfrac{a_s^{(k)}}{\|a_{s}\|}\frac{\lambda_s^{k,l}}{\Lambda^k_s}\left( c - \sum_{h=0}^{c} (c-h)  \frac{(\Lambda_s^k)^{h}}{h!}e^{-\Lambda_s^k} \right) - \dfrac{a_s^{(l)}}{\|a_s\|};\\
		f_{2l}(x_s) &= \sum_{k=1}^{m} \dfrac{a_s^{(k)}}{\|a_{s}\|} \mu_s^{k,l} -  \sum_{k=1}^{m} \dfrac{a_s^{(k)}}{\|a_{s}\|} \frac{\lambda_s^{k,l}}{\Lambda^k_s}\left( c - \sum_{h=0}^{c} (c-h)  \frac{(\Lambda_s^k)^{h}}{h!}e^{-\Lambda_s^k} \right);\\
		f_{3l}(x_s) &= \dfrac{a_s^{(l)}}{\|a_s\|};
		\end{align}	with
		\begin{align}
			\lambda_s^{k,l}&:= \lambda_{kl}\left(\pi_l-a^{(l)}_s - u^{(l)}_s\right); \quad \Lambda^k_s:= \sum_{l=1}^{m}\lambda^{k,l}_s \quad
		\text{and} \quad \mu^{k,l}_s:= \lambda_{kl}(\pi_l - a^{(l)}_s - b^{(l)}_s - u^{(l)}_s).
		\end{align}
		For $s \in [t_0,1],  f(x_s) = f(x_{t_0}).$		
	\end{thrm}

\begin{rmrk}Notice that in this model, the time corresponds to the fraction of the population interviewed. The time $t_0$ is the first time at which $|a_t|$ reaches $0$ and can be seen as the proportion of the population interviewed when there is no more coupon to keep the CRS going. Necessarily, $t_0\leq 1$. We see that $\|a_t\|=0$ only if $a_t^{(1)} = \ldots = a_t^{(m)}=0$. It implies that $f(x_t) = 0, \forall t \in [t_0,1]$. Then, the solution of the system of ODEs \eqref{th:main} becomes constant over the interval $[t_0, 1]$.
	\end{rmrk}
	
	The rest of this paper is organized in the following manner. First, in Section 2, we give a precise description of the chain-referral process on a SBM random graph. This relies heavily on the structure of the random graph that we construct progressively when the exploration process spreads on it. In {Section 3}, we prove the limit theorem. The proof uses limit theory of c\`adl\`ag semi-martingale vector processes {equipped with} Skorokhod topology (see \cite{ethierkurtz}) and Poisson approximations (see \cite{barbour1992}). Then in {Section 4}, {we present simulation results} of the stochastic process and the solution of the system of limiting ODEs. We conclude with some discussions on the impacts of changing parameters of the models on the evolution of the chain-referral process.
	
	\section{Definition of the chain-referral process}
	
	Let us describe the dynamics of $X=(X_n)_{n\in\N}$. Recall that $\|A_n\|:=\sum_{l=1}^{m} A_n^{(l)}$ is the {total number of} individuals having coupons but {who have not yet been interviewed}. We start with $A_0$ seeds, whose types are chosen independently according to $\pi$. {$A_0$ is an m-dimensional random vector with multinomial distribution $\mathcal{M}(\lfloor \|a_0\| N \rfloor;\pi_1,...,\pi_m)$, \ie $\P\big((A_0^{(1)},...,A_0^{(m)})=(k_1,...,k_m)\big) = \pi_1^{k_1}...\pi_m^{k_m}, k_i \in \N$ such that $\sum_{i=1}^m k_i = \lfloor \|a_0\| N \rfloor$ and Assumption \ref{assp:initialcondition} is satisfied.} Also $B_0=U_0=(0,...,0)$ and we set $X_0=(A_0,B_0,U_0)$.
	
	{We now define $X_n$ given the state $X_{n-1}$ previous to the $n^{th}$-interview and given the number $N_1, \ldots, N_m$ of nodes of each type.} At step $n \geq 1$, after the $n^{th}$-interview, the type of the upcoming interviewee is chosen uniformly at random according to the number of active coupons of each type in the present time. To choose the type of the next interviewee, we define an $m$-dimensional vector $I_n:=(I_n^{(1)},...,I_n^{(m)})$, which takes value $1$ at coordinate $l$ and $0$ elsewhere if the $n^{\text{th}}$ interviewee belongs to block $l$. This $n^{th}$-interviewee is chosen uniformly among the $\|A_{n-1}\|$ active coupons of $m$ types \ie $I_n$ has multinomial distribution
	\begin{equation} \label{eq:defofIn}
	I_n = (I_n^{(1)},...,I_n^{(m)}) \stackrel{(d)}{=} \mathcal{M}\left(1;\frac{A_{n-1}^{(1)}}{\|A_{n-1}\|},...,\frac{A_{n-1}^{(m)}}{\|A_{n-1}\|}\right).
	\end{equation}
	If the chosen one belongs to block $[l]$, $A_n^{(l)}$ is reduced by 1 and a number of new coupons distributed are added up, depending on how many new contacts he/she has. {In the meantime, the number of interviewees of type $l$ is increased by $1$. \ie $U_n^{(l)} = U_{n-1}^{(l)} + I_n^{(l)}$.} Among the new contacts of the $n^{\text{th}}-$interviewee, define $H_n^{(l)}$ the number of new contacts of type $l$, who have not been mentioned before; $K_n^{(l)}$ the number of new contacts of type $l$ whose identities are already known but who are still inactive. The $H^{(l)}_n$ new connections are chosen independently among $N_l- A^{(l)}_{n-1}- B_{n-1}^{(l)} - U_n^{(l)}$ individuals in the hidden population {where probability of each successful connection is} $\sum_{k=1}^{m}I_n^{(k)}p_{kl}$. {Hence, conditioning on $(N_1,\ldots, N_m), X_{n-1}$, the random variable $H^{(l)}_n$ follows the binomial distribution:}
	\begin{align} \label{distributionofHn}
	H_n^{(l)} &\stackrel{(d)}{=} \mathcal{B}\textit{in}\left(N_l-A_{n-1}^{(l)}-B_{n-1}^{(l)}-U_n^{(l)},\sum_{k=1}^{m}I_n^{(k)}p_{kl}\right).
	\end{align}
	And the $K^{(l)}_n$ individuals are chosen independently of $H_n^{(l)}$ from $B_{n-1}^{(l)}$ individuals and independently of the others with probability $\sum_{k=1}^{m}I_n^{(k)}p_{kl}$. In that way, {conditioning on $(N_1, \ldots, N_m), X_{n-1}$,  $K_n^{(l)}$ also has the binomial distribution:}
	\begin{equation}
	K_n^{(l)} \stackrel{(d)}{=} \mathcal{B}\textit{in} \left(B_{n-1}^{(l)},\sum_{k=1}^{m}I_n^{(k)}p_{kl}\right).
	\end{equation}
	In total, there are $Z_n:=H_n+K_n$ candidates, who can possibly receive coupons at step $n$. Notice that, {conditioning on $(N_1, \ldots, N_m), X_{n-1}$, $(H^{(l)}_n)_{l=1,...,m}$ and $(K^{(l)}_n)_{l=1,...,m}$ are independent,} henceforth,
	\begin{equation}
	Z_n^{(l)}\stackrel{(d)}{=} \mathcal{B}\textit{in}\left(N_l-A_{n-1}^{(l)}-U_n^{(l)},\sum_{k=1}^{m}I_n^{(k)}p_{kl}\right).
	\end{equation}
	Let $C_n=(C_n^{(1)},...,C_n^{(m)})$ ($l=1,...,m$) be the numbers of coupons that are distributed at step $n$. By the setting of the survey, the total coupons $|C_n|$ must be maximum $c$. If {the number $Z_n$ of candidates} is less than or equal to $c$, we deliver exactly $Z_n$ coupons. Otherwise, we choose new people to be enrolled in the study by an $m-$dimensional random variable $C'^{(l)}_n= (C'^{(1)}_n, ...,C'^{(m)}_n)$ having the multivariate hypergeometric distribution with parameters $(m; c, (Z^{(1)}_n,...,Z^{(m)}_n))$ and the support $\{(c_1,...,c_m) \in \N^m: \forall l \leq m, c_l \leq Z^{(l)}_n, \sum\limits_{l=1}^{m}c_i = c\}$, that is
	\begin{equation*}
	\P \left((C'^{(1)}_n,...,C'^{(m)}_n)=(c_1,...,c_m)\right)= \dfrac{\prod\limits_{l=1}^{m} \binom{Z_n^{(l)}}{c_l}}{\binom{\sum_{l=1}^{m}Z_n^{(l)}}{c}}.
	\end{equation*}
	In another words,
	\begin{equation}
	C^{(l)}_n:= \begin{cases}
	Z^{(l)}_n & \quad \text{if} \quad \sum_{l=1}^{m}Z_n^{(l)} \leq c\\
	C'^{(l)}_n& \quad \text{otherwise}
	\end{cases}.
	\end{equation}
Let define by
	\begin{equation}
	    n_0:= \inf\{ n \in \{1,...,N\}, A_n=0\}
	\end{equation}
	the first step that $|A_n|$ reaches zero.
	{The dynamics of $X_n$ can be described by the following recursion:}
	\begin{align}
	&\begin{cases}\label{def:ABU}
	A_n &= A_{n-1} - I_n + C_n\\
	B_n &= B_{n-1} + H_n - C_n\\
	U_n &= \sum\limits_{i=1}^{n} I_i
	\end{cases}, \quad \text{for} \quad n \in \{1,...,n_0\}\\
	\text{and} \quad & \quad X_n \quad =  X_{n-1} \quad \text{when} \quad n> n_0. \nonumber
	\end{align}
The random network is progressively discovered when the referrals chain process explores it.
	\begin{prpstn}
		Consider the discrete-time process $(X_n)_{1\leq n \leq N}$ defined in \eqref{def:ABU}. For $n \in \mathbb{N}$, we denote by $\mathcal{F}_{n}:= \sigma\big(\{ X_i, i \leq n , (N_1, \ldots, N_m)\}\big)$ the canonical filtration associated with $(X_n)_{1\leq n \leq N}$. Then the process $(X_n)_{n}$ is an inhomogeneous Markov chain with respect to the filtration $(\mathcal{F}_n)_n$.
	\end{prpstn}
	\begin{proof}
		The proposition is deduced from the recursion \eqref{def:ABU} of $(X_n)_{1\leq n\leq N}$ and the fact {that the random variables $C_n, I_n, H_n$ are defined conditionally on $X_{n-1}$ and $(N_1, \ldots, N_m)$}. {The fact that the Markov process is inhomogeneous comes from the setting of the CRS survey: there is no replacement in the recruitment procedure. For example, when $m=1$, the definition of $H^{(l)}_n$ in \eqref{distributionofHn} depends on time as $U_n^{(l)}=n$.}
	\end{proof}
	
	\section{Asymptotic behavior of the chain-referral process}
	Let us now consider the renormalized chain-referral process given in \eqref{eq:renorm} in the time interval $[0,t_0]$. The main theorem (Theorem \ref{th:main}) shows the convergence of the sequence $(X_\cdot^N)_{N}$ to a deterministic process. For this, we look for an expression of the equations \eqref{def:ABU} as a vector of semi-martingales. We start by writing the Markov chain $(X_n)_{1 \leq n \leq N}$ as the sum of its increments in discrete time.
	\begin{equation*}
	X_n= X_0 + \sum_{i=1}^{n} (X_i-X_{i-1})= \begin{pmatrix}
	A_0\\B_0\\U_0
	\end{pmatrix} +  \sum_{i=1}^{n}\begin{pmatrix}
	C_i-I_i\\H_i-C_i\\I_i
	\end{pmatrix}.
	\end{equation*}
	{Each element of the increment $X_{n+1} - X_{n}$ are binomial variables conditioned on all the events having been occurring until step $n$. When we fix $n$ and let $N$ tend to infinity, the conditional binomial random variables can be approximated by some Poisson random variables. The normalization $X^N_t$ of $X_n$ becomes:}
	\begin{align*}
		X^N_t = \frac{1}{N} \begin{pmatrix}
		A_0\\B_0\\U_0
		\end{pmatrix}  + \frac{1}{N}\sum_{i=1}^{\lfloor Nt \rfloor}\begin{pmatrix}
		C_i-I_i\\H_i-C_i\\I_i
		\end{pmatrix}.
	\end{align*}
{The Doob decomposition of the renormalized processes $(X_t^N)_{t\in[0,t_0]}$ given in Section \ref{sec:Doob} consists of a finite variation process and an $\mathbb{L}^2$-martingale. We use Aldous criteria (conditionally on the past see e.g. \cite{ethierkurtz,Metivier}) to show the tightness of the distributions of these processes in Section \ref{sec:tightness}. Once the tightness is established, we identify the limiting values of this tight sequence and finally we prove that the limiting values of all converging subsequences are the same, hence it is the limit of processes $(X_\cdot^N)_N$.} This proves Theorem \ref{th:main}. \\

	Denote by $(\mathcal{F}^N_t)_{t\in [0,1]}:=(\mathcal{F}_{\lfloor Nt \rfloor})_{t\in [0,1]}$ the canonical filtration associated to $(X^N_t)_{t\in [0,1]}$.


	\subsection{Doob's decomposition}\label{sec:Doob}
	\begin{lmm}
		The process $(X_t^N)_{t\in[0,1]}$ admits the Doob's decomposition: $X^N_t=X_0^N + \Delta^N_t + M^N_t$, $X_0^N = \frac{1}{N}X_0$. $(\Delta_t^N)_{t \in [0,1]}$ is an $\mathcal{F}_t^N-$predictable process defined by
		\begin{equation}\label{eq:Delta}
		\Delta_t^N = \begin{pmatrix}
		\Delta^{N,1}_t\\
		\Delta^{N,2}_t\\
		\Delta^{N,3}_t
		\end{pmatrix}
		=\frac{1}{N}\sum_{n=1}^{\lfloor Nt \rfloor}
		\begin{pmatrix}
		\E[C_n-I_n|\mathcal{F}_{n-1}]\\
		\mathbb{E}[H_n-C_n|\mathcal{F}_{n-1}]\\
		\mathbb{E}[I_n|\mathcal{F}_{n-1}]
		\end{pmatrix};
		\end{equation}
		$(M_t^N)_{t\in [0,1]}$ is an $\mathcal{F}_t^N-$ square integrable centered martingale with quadratic variation process $(\langle M_\cdot^N \rangle_t)_{t\in [0, 1]}$ given by: for every $(l,k) \in \{1,...,m\}^2$,
		\begin{equation}\label{bracketofmartingale}
		\langle M_\cdot^{(l),N},M_\cdot^{(k),N}\rangle_t = \frac{1}{N^2}\sum_{n=1}^{\lfloor Nt \rfloor} \E\left[\left( X^{(l)}_n - \E[X^{(l)}_n |\mathcal{F}_{n-1}] \right)\left( X^{(k)}_n - \E[X^{(k)}_n|\mathcal{F}_{n-1}] \right)^{T} \bigg|\mathcal{F}_{n-1}\right], \quad t\in [0,1]
		\end{equation}
		where $X$ is a column vector and $X^T$ is its transpose.
	\end{lmm}
	
	\begin{proof}
		In order to obtain the Doob's decomposition, we write for $t\in [0,1]$,
		\begin{align*}
		X_t^N& =\frac{X_0}{N} + \frac{1}{N}\sum_{n=1}^{\lfloor Nt \rfloor} (X_n - X_{n-1}) \\& = X_0^N + \frac{1}{N}\sum_{n=1}^{\lfloor Nt \rfloor} \E[X_n-X_{n-1}| \mathcal{F}_{n-1}] + \frac{1}{N}\sum_{n=1}^{\lfloor Nt \rfloor} (X_n - X_{n-1} - \E[X_n-X_{n-1} | \mathcal{F}_{n-1}])\\
		& = X_0^N + \Delta_t^N + M_t^N.
		\end{align*}
		It is clear that the conditional expectations above are all well-defined since the components of $X_n$ and $X_{n-1}$ are all bounded by $N$, that $\Delta^N_t$ is $\mathcal{F}^N_{t}-$predictable and that $(M^N_t)_{t\in [0,1]}$ is an $\mathcal{F}_{t}^N-$martingale. We first check that $(\Delta_\cdot^N)_N$ is a sequence of finite variation processes and then we can conclude  that $X^N_t= X_0^N + \Delta_t^N + M^N_t$ is the Doob's decomposition.\\
		Denote by $\lambda:= \displaystyle\max_{l,k \in \{1,...,m\}}\lambda_{kl}$.
		Notice that
		\begin{align}
		\|\E[A_n-A_{n-1}|\mathcal{F}_{n-1}\|&= \|\E[C_n-I_{n}|\mathcal{F}_{n-1}\|\leq c \label{incrementsofXn1},\\
		\|\E[B_n-B_{n-1}|\mathcal{F}_{n-1}\|&= \|\E[H_n-C_{n}|\mathcal{F}_{n-1}\| \leq m(\displaystyle\max_{l,k \in \{1,...,m\}}\lambda_{kl})+c= m\lambda + c \label{incrementsofXn2},\\
		\|\E[U_n-U_{n-1}|\mathcal{F}_{n-1}\| &\leq 1, \label{incrementsofXn3}
		\end{align}
		then $\|\E[X_n-X_{n-1}|\mathcal{F}_{n-1}]\| \leq 2c+ m \lambda +1$. So the total variation of $(\Delta_t^N)_{t\in[0,1]}$ is
		\begin{equation*}
		V^N(\Delta_t^N) = \frac{1}{N}\sum_{n=1}^{\lfloor Nt \rfloor}  \|\Delta^N_{nt/N}-\Delta^N_{(n-1)t/N}\| = \frac{1}{N}\sum_{n=1}^{\lfloor Nt \rfloor} \|\E[X_n-X_{n-1}| \mathcal{F}_{n-1}]\| \leq (2c+ m\lambda +1)t,
		\end{equation*}
		which is finite. It follows that $(\Delta_t^N)_{t\in[0,1]}$ is an $\mathcal{F}^N_t-$ predictable with finite variations.
		
		The quadratic variation of $(M_t^N)_{t\in [0,1]}$ is computed as follow. For every $k,l=1,...,m$
		\begin{align*}
		M^{(l),N}_t\left(M^{(k),N}_t\right)^T = \frac{1}{N^2}\sum_{n=1}^{\lfloor Nt \rfloor} \left( X^{(l)}_n- X^{(l)}_{n-1} - \E[X^{(l)}_n- X^{(l)}_{n-1}|\mathcal{F}_{n-1}] \right)\left( X^{(k)}_n- X^{(k)}_{n-1} - \E[X^{(k)}_n- X^{(k)}_{n-1}|\mathcal{F}_{n-1}] \right)^{T} \\
		+ \frac{1}{N^2}\sum_{n=1}^{\lfloor Nt \rfloor}\sum_{\substack{n'=1\\ n'\neq n}}^{\lfloor Nt \rfloor} \left( X^{(l)}_n- X^{(l)}_{n-1} - \E[X^{(l)}_n- X^{(l)}_{n-1}|\mathcal{F}_{n-1}] \right)\left( X^{(k)}_{n'}- X^{(k)}_{n'-1} - \E[X^{(k)}_{n'}- X^{(k)}_{n'-1}|\mathcal{F}_{n'-1}] \right)^{T}\\
		=:L^N_t + L'^N_t. \hspace*{33em}
		\end{align*}
		The term $L'^N_t$ is an $\mathcal{F}^N_{t}-$martingale since whenever $n'<n$, $\left( X^{(k)}_{n'}- X^{(k)}_{n'-1} - \E[X^{(k)}_{n'}- X^{(k)}_{n'-1}|\mathcal{F}_{n'-1}] \right) $ is $\mathcal{F}_{n-1}-$measurable. To see that the quadratic variation of $M^N_t$ has the form \eqref{bracketofmartingale}, we write the term $L^N_t$ as follows:
		\begin{align*}
		L^N_t &:= \frac{1}{N^2}\sum_{n=1}^{\lfloor Nt \rfloor} \E\left[ \left( X^{(l)}_n - \E[X^{(l)}_n |\mathcal{F}_{n-1}] \right)\left( X^{(k)}_n - \E[X^{(k)}_n|\mathcal{F}_{n-1}] \right)^{T} \bigg|\mathcal{F}_{n-1} \right] \\
		&\quad+ \frac{1}{N^2}\sum_{n=1}^{\lfloor Nt \rfloor} \left( X^{(l)}_n - \E[X^{(l)}_n |\mathcal{F}_{n-1}] \right)\left( X^{(k)}_n - \E[X^{(k)}_n|\mathcal{F}_{n-1}] \right)^{T} \\
		& \quad - \frac{1}{N^2}\sum_{n=1}^{\lfloor Nt \rfloor} \E\left[\left( X^{(l)}_n - \E[X^{(l)}_n |\mathcal{F}_{n-1}] \right)\left( X^{(k)}_n - \E[X^{(k)}_n|\mathcal{F}_{n-1}] \right)^{T}\bigg|\mathcal{F}_{n-1}\right]\\
		&= \frac{1}{N^2}\sum_{n=1}^{\lfloor Nt \rfloor} \E\left[\left( X^{(l)}_n - \E[X^{(l)}_n |\mathcal{F}_{n-1}] \right)\left( X^{(k)}_n - \E[X^{(k)}_n|\mathcal{F}_{n-1}] \right)^{T}  \bigg|\mathcal{F}_{n-1}\right] + L''^N_t = \langle M^N\rangle_t + L''^N_t.
		\end{align*}
		As a result,
		\begin{equation}
		M^{(l),N}_t\left(M^{(k),N}_t\right)^T= \langle M^N\rangle_t + L'^N_t + L''^N_t.
		\end{equation}
		Because both $L'^N_t$ and $L''^N_t$ are $\mathcal{F}^N_t-$martingale, $L'^N_t + L''^N_t$ is an $\mathcal{F}^N_t-$martingale as well. The term $ (\langle M^N\rangle_t)_t$ is $\mathcal{F}^N_t-$adapted with the variation
		\begin{align}\label{eq:variance}
		V^N(\langle M^N_\cdot \rangle_t) &= \frac{1}{N^2} \sum_{n=1}^{\lfloor Nt \rfloor}\sum_{k,l=1}^{m}\left\| \E\left[\left( X^{(l)}_n - \E[X^{(l)}_n |\mathcal{F}_{n-1}] \right)\left( X^{(k)}_n - \E[X^{(k)}_n|\mathcal{F}_{n-1}] \right)^{T} \bigg|\mathcal{F}_{n-1}\right]\right\|.
\end{align}The integrand in the right hand side is the conditional covariance between $X_n^{(l)}$ and $X_n^{(k)}$ conditionally to $\mathcal{F}_{n-1}$. Because $X_n^{(l)}$ and $X_n^{(k)}$ are vectors, this covariance is a matrix of size $3\times 3$ and for $1\leq i,j\leq 3$, the term $(i,j)$ of this matrix is:
\begin{multline*}\E\left[\left( X^{i,(l)}_n - \E[X^{i,(l)}_n |\mathcal{F}_{n-1}] \right)\left( X^{j,(k)}_n - \E[X^{j,(k)}_n|\mathcal{F}_{n-1}] \right)\bigg|\mathcal{F}_{n-1}\right] \\ \leq   \left(\Var(X_n^{i,(l)} -X_{n-1}^{i,(l)}|\mathcal{F}_{n-1}) \right)^{1/2} \left( \Var(X_n^{j,(k)} -X_{n-1}^{j,(k)}|\mathcal{F}_{n-1}) \right)^{1/2},\end{multline*}
by the Cauchy-Schwarz inequality. Thus:
\begin{align*}
V^N(\langle M^N_\cdot \rangle_t)		&\leq \frac{1}{N^2} \sum_{n= 1}^{\lfloor Nt \rfloor} \sum_{k,l=1}^{m}\left| \sum_{i,j=1}^{3} \left(\Var(X_n^{i,(l)} -X_{n-1}^{i,(l)}|\mathcal{F}_{n-1})\right)^{1/2} \left(\Var(X_n^{j,(k)} -X_{n-1}^{j,(k)}|\mathcal{F}_{n-1})\right)^{1/2} \right|,
		\end{align*}
		where $(X_n^{1,(l)},X_n^{2,(l)},X_n^{3,(l)})=(A^{(l)}_n,B_n^{(l)},U_n^{(l)})$. By Cauchy-Schwarz's inequality, we have
		\begin{align}
			\sum_{i,j=1}^{3} \left(\Var(X_n^{i,(l)} -X_{n-1}^{i,(l)}|\mathcal{F}_{n-1})\right)^{1/2} \left( \Var(X_n^{j,(k)} -X_{n-1}^{j,(k)}|\mathcal{F}_{n-1}) \right)^{1/2}  \hspace*{2.5in} \nonumber\\
			= \left(\sum_{i=1}^{3} \left(\Var(X_n^{i,(l)} -X_{n-1}^{i,(l)}|\mathcal{F}_{n-1}) \right)^{1/2} \right) \left( \sum_{j=1}^{3} \left(\Var(X_n^{j,(k)} -X_{n-1}^{j,(k)}|\mathcal{F}_{n-1})\right)^{1/2}\right) \nonumber\\
		  \leq \frac{3}{2} \sum_{i =1}^{3 }\left( \Var(X_n^{i,(l)} -X_{n-1}^{i,(l)}|\mathcal{F}_{n-1}) + \Var(X_n^{i,(k)} -X_{n-1}^{i,(k)}|\mathcal{F}_{n-1}) \right). \hspace*{0.85in} \label{ineq:variation}
		\end{align}
		 From \eqref{incrementsofXn1}-\eqref{incrementsofXn3} and by Cauchy-Schwarz's inequality, we obtain the following inequalities
		\begin{equation}\label{variationofbracket}
		\Var(C_n^{(l)} -I_n^{(l)}|\mathcal{F}_{n-1}) \leq c^2, \quad \Var(H_n^{(l)} -C_{n}^{(l)}|\mathcal{F}_{n-1}) \leq 2(\displaystyle\max_{l,k \in \{1,...,m\}} \lambda_{lk}^2 + c^2), \quad \Var(I_n^{(l)}|\mathcal{F}_{n-1}) \leq 1.
		\end{equation}
		As a consequence,
		\begin{align*}
		V^N\big(\langle M^N_\cdot \rangle_t \big) \leq \frac{1}{N^2} \sum_{n= 1}^{\lfloor Nt \rfloor} 3m^2(c^2+ 2(\displaystyle\max_{l,k \in \{1,...,m\}} \lambda_{lk}^2 + c^2) + 1) \leq \frac{1}{N}3m^2 (3c^2 + 2\lambda^2+1)< \infty.
		\end{align*}
		Thus, the proof of the Lemma is completed.
	\end{proof}
	
	\subsection{Tightness of the renormalized process}\label{sec:tightness}
	\begin{lmm}\label{lem:tightness}
		The sequence of processes $(X_\cdot^N)_N$ is tight in the Skorokhod space $\mathcal{D}([0,1],[0,1]^{3\times m})$.
	\end{lmm}
	\begin{proof}
		
		To prove the tightness of $(X_\cdot^N)_N$, we use the criteria of tightness {for semi-martingales in \cite[Theorem 2.3.2 (Rebolledo)]{Metivier}: first, we verify the marginal tightness of each sequence $(X^N_t)_N$ for each $t \in [0,1]$, then we show the tightness for each process in the Doob's decomposition of $X^N$, the finite variation process $(\Delta^N)_N$ and the quadratic variation of the martingale $(M^N)_N$}. {For any $t\in [0,1]$, the tightness of marginal sequence $(X_t^N)_N$ is easily deduced from the compactness of a sequence of random variables taking values in a compact set $[0,1]^{3\times m}$.} 		
		{Since the sequence of martigales $(M^N)_N$ is proved to be convergent (to zero) in $\mathbb{L}^2$ as $N \rightarrow \infty$ (which is done by Proposition \ref{prop:cvofmartingale}), we have the tightness of $(M^N)_N$.} Thus, it is sufficient to check the tightness condition for the modulus of continuity of $(\Delta^N)_N$ (see, \eg, \cite[Theorem 13.2, p.139]{billingsley}).\\
		For all $0<\delta<1$ and for every $s,t \in [0,1]$ such that $|t-s|<\delta$, we have that
		\begin{align*}
		\|\Delta^N_t - \Delta^N_s\| = \left\|\frac{1}{N} \sum_{n=\lfloor Ns \rfloor + 1}^{\lfloor Nt \rfloor} \E[X_n- X_{n-1}|\mathcal{F}_{n-1}]\right\|& \leq \frac{1}{N} \sum_{n=\lfloor Ns \rfloor + 1}^{\lfloor Nt \rfloor} \|\E[X_n- X_{n-1}|\mathcal{F}_{n-1}]\|.
		\end{align*}
		By \eqref{incrementsofXn1}-\eqref{incrementsofXn3}, we get
		\begin{equation*}
		\|\Delta^N_t - \Delta^N_s\| \leq \frac{\lfloor Nt \rfloor - \lfloor Ns \rfloor}{N}(c+ m\lambda +c +1) \leq (2c+ m\lambda +1) \left(\delta+\frac{1}{N}\right).
		\end{equation*}
		Thus, for each $\varepsilon>0$, choose $\delta_0\leq  \frac{\varepsilon}{2(2c+m\lambda + 1)}$, we have that
		\begin{align*}
		\P\left(\sup_{\substack{|t-s|<\delta\\  0\leq s<t \leq 1}}\|\Delta^N_t-\Delta^N_s\| >\varepsilon\right) =0, \quad \forall \delta \leq \delta_0, \forall N > \frac{1}{\delta_0},
		\end{align*}
		which allows us to conclude that the sequence $(\Delta_\cdot^N)_N$ is tight and finishes the proof of the lemma.
	\end{proof}

To complete the proof of Lemma \ref{lem:tightness}, we now prove that:

	\begin{prpstn} \label{prop:cvofmartingale}
	The sequence of martingale $(M_\cdot^N)_N$ converges to $0$ in $\mathbb{L}^2$ as $N$ goes to infinity.
\end{prpstn}

\begin{proof}		Consider the quadratic variation of $(M_\cdot^N)_N$:
	According to the fomula \eqref{bracketofmartingale}, we apply the Cauchy-Schwarz's inequlity and then use the inquality \eqref{ineq:variation} to obtain that for every $t \in [0,1]$,
	\begin{align*}
	\|\langle M^{(l),N},M^{(k),N}\rangle_t \|
	&=\left\|\frac{1}{N^2} \sum_{n= 1}^{\lfloor Nt \rfloor} \E\left[ \left( X^{(l)}_n - \E[X^{(l)}_n |\mathcal{F}_{n-1}] \right)\left( X^{(k)}_n - \E[X^{(k)}_n|\mathcal{F}_{n-1}] \right)^{T} \bigg|\mathcal{F}_{n-1} \right]\right\| \\
	&\leq \frac{1}{N^2} \sum_{n= 1}^{\lfloor Nt \rfloor} \left| \sum_{i,j=1}^{3} \left(\Var(X_n^{i,(l)} -X_{n-1}^{i,(l)}|\mathcal{F}_{n-1})\right)^{1/2} \left( \Var(X_n^{j,(k)} -X_{n-1}^{j,(k)}|\mathcal{F}_{n-1}) \right)^{1/2} \right| \\
	& \leq \frac{1}{N^2} \sum_{n= 1}^{\lfloor Nt \rfloor} \frac{3}{2} \sum_{i =1}^{3 }\left( \Var(X_n^{i,(l)} -X_{n-1}^{i,(l)}|\mathcal{F}_{n-1}) + \Var(X_n^{i,(k)} -X_{n-1}^{i,(k)}|\mathcal{F}_{n-1}) \right),
	\end{align*}
	where $(X_n^{1,(l)},X_n^{2,(l)},X_n^{3,(l)})=(A^{(l)}_n,B_n^{(l)},U_n^{(l)})$.  From \eqref{incrementsofXn1}-\eqref{incrementsofXn3} and \eqref{variationofbracket}, we deduce that
	\begin{equation}
	\|\langle M^N_\cdot\rangle_t \| \leq \frac{1}{N^2} \sum_{n=1}^{\lfloor Nt \rfloor} \frac{3m^2}{2}\left( c^2+ 2(\displaystyle\max_{l,k \in \{1,...,m\}} \lambda_{lk}^2 + c^2) + 1 \right) \leq \frac{1}{N}\frac{3m^2}{2} (3c^2 + 2\lambda^2+1)t.
	\end{equation}
	Applying the Doob's inequality for martingale, for every $t\in [0,1]$, we have
	\begin{equation*}
	\E\bigg[\max_{0\leq s \leq t}\|M^N_s\|^2\bigg] \leq 4\E\bigg[\|\langle M^N_\cdot\rangle_t\|\bigg] \leq \frac{1}{N}6m^2 (3c^2 + 2\lambda^2+1) \rightarrow 0 \quad \text{as } \quad N \rightarrow \infty.
	\end{equation*}This concludes the proof of Prop. \ref{prop:cvofmartingale} and hence of Lemma \ref{lem:tightness}.
\end{proof}
	
	\subsection{Identify the limiting value}\label{sec:identif}
	Since the sequence $(X_\cdot^N)_N$ is tight, for any limiting value $x=(a,b,u)$ of the sequence $(X^N)_N$, there exists an increasing sequence $(\varphi_N)_N$ in $\N$ such that $(X_\cdot^{\varphi_N})_N$ converges in distribution to $x$ in $\mathcal{D}([0,1],[0,1]^{3\times m})$. Because the sizes of the jumps converge to zero with $N$, the limit is in fact in $\mathcal{C}([0,1],[0,1]^{3\times m})$. We want to identify that limit. In order to simplify the notations, we also write the subsequence $(X_\cdot^{\varphi_N})_N$ as $(X_\cdot^N)_N= (A_\cdot^N,B_\cdot^N,U_\cdot^N)_N$.\\
	
	We consider separately the martingale and finite variation parts. Proposition \ref{prop:cvofmartingale} implies that the sequence martingale $(M_\cdot^N)_N$ converges to $0$ in distribution {and hence $(M^N)_N$ converges to zero in probability.} It remains to find the limit of the finite variation process $(\Delta_\cdot^N)_N$ given in Equation \eqref{eq:Delta} {and prove that the limit found is the same (which is done later in the proof for the uniqueness of the system of the ODEs \eqref{th:main}) for every convergent subsequence extracted from the tight sequence $(X^N)_N$.}
	
	\begin{prpstn}\label{prop:conv-Delta}When $N$ goes to infinity, we have the following convergences in distribution in $\mathcal{D}([0,1],[0,1]^{3\times m})$:
		\begin{align}
		&	\frac{1}{N}\sum_{n=1}^{\lfloor Nt \rfloor} \E[C^{(l)}_n|\mathcal{F}_{n-1}] \stackrel{(d)} {\rightarrow}\int_{0}^{t}\left\{ \sum_{k=1}^{m} \dfrac{a_s^{(k)}}{\|a_s\|}\frac{\lambda_s^{k,l}}{\Lambda^k_s}\left( c - \sum_{h=0}^{c} (c-h)  \frac{(\Lambda_s^k)^{h}}{h!}e^{-\Lambda_s^k} \right) \right\}ds, \label{limitofC_n}\\
		&	\frac{1}{N}\sum_{n=1}^{\lfloor Nt \rfloor} \E[H^{(l)}_n|\mathcal{F}_{n-1}] \stackrel{(d)} {\longrightarrow}\int_{0}^{t}\sum_{k=1}^{m} \dfrac{a_s^{(k)}}{\|a_s\|}\mu^{k,l}_s ds, \label{lim:limitofHn}\\
		& 	\frac{1}{N} \sum_{n= 1}^{\lfloor Nt \rfloor}\E[I^{(l)}_n|\mathcal{F}_{n-1}]=\frac{1}{N} \sum_{n= 1}^{\lfloor Nt \rfloor}\left(\frac{A_{n-1}^{(l)}}{N}\right)\bigg/\left(\frac{\|A_{n-1}\|}{N}\right) \stackrel{(d)} {\longrightarrow} \displaystyle\int\limits_{0}^{t} \frac{a_s^{(l)}}{\|a_s\|}ds, \label{lim:limitofIn}
		\end{align} {where $\lambda_s^{k,l}, \Lambda_s^{k}, \mu_s^{k,l}$ are defined as in Theorem \ref{th:main}.}
		This provides the convergence of $(\Delta_\cdot^N)_N$ to a solution $x_.$ of \eqref{limiting_eq}.
		
	\end{prpstn}
Since the limits are deterministic, the convergences hold in probability. Moreover the uniqueness of the solution of \eqref{limiting_eq} will be proved later, which will imply the convergence of the whole sequence $(X_\cdot^N)_N$ to this solution.\\
	
\begin{proof} Recall that since the sequence $(X^N_.)_{N}$ is tight, we have extracted a converging subsequence also denoted by $(X^N_.)_N$ of which we study the limit.\\

	The proof of the Proposition \ref{prop:conv-Delta} is separated into three steps. \\
	\textbf{Step 1:} We consider the most complicated term $\E[C_n|\mathbb{F}_{n-1}]$. We prove that:
	for each $l\in \{0,...,m\}$,
		\begin{equation}\label{EC_n_in_lemma}
		\left| \E[C_n^{(l)}|\mathcal{F}_{n-1}] - \frac{\lambda_n^{(l)}}{\Lambda_n}\left( c - \sum_{h=0}^{c} (c-h)  \frac{(\Lambda_n)^{h}}{h!}e^{-\Lambda_n} \right) \right| \leq \frac{m(c+1)\lambda}{N},
		\end{equation}
		where
		\begin{equation}
		\lambda_n^{(l)}:=\left(\sum_{k=1}^{m}I_n^{(k)}\lambda_{kl}\right)\left( \frac{N_l}{N}-\frac{A_{n-1}^{(l)}}{N}-\frac{U_n^{(l)}}{N}\right) \quad  \quad \text{and} \quad \quad \Lambda_n:= \sum_{j=1}^{m}\lambda_n^{(j)}.
		\end{equation}\\
		
		{Notice that $\Lambda_n = 0$ only if for each $l \in \{1, \ldots, m\}$, $\lambda_n^{(l)} = 0$. It happens when $A_{n-1}^{(l)} +U_{n}^{(l)} = N_l$, meaning that all the nodes of type $l$ have been discovered. In this case, $C_n^{(l)} =0$ and \eqref{EC_n_in_lemma} is satisfied.}\\
	    Let us write
		\begin{align}\label{EC_n}
		\mathbb{E}[C_n^{(l)}|\mathcal{F}_{n-1}]= \mathbb{E}\left[Z^{(l)}_n \mathbb{1}_{\sum_{j=1}^{m}Z^{(j)}_n\leq c}\big|\mathcal{F}_{n-1}\right] + \mathbb{E}\left[\frac{cZ^{(l)}_n}{\sum_{j=1}^{m}Z^{(j)}_n}\mathbb{1}_{\sum_{j=1}^{m}Z^{(j)}_n>c}\bigg|\mathcal{F}_{n-1}\right].
		\end{align}
		For every $l=1,...,m$ and every fixed $n$, {when all the parameters are positive}, we have that $(N_l-A_{n-1}^{(l)}-U_n^{(l)})\underset{\text{a.s.}}{\overset{N \rightarrow \infty}{\longrightarrow}} +\infty$. Then we work conditionally on $N_l, A_{n-1}^{(l)}, U_n^{(l)}$ and $I_n^{(l)}$ and use {the Poisson approximation (\eg see Equation (1.23) and Theorem 2.A, 2.B by Barbour, Holst and Janson in \cite{barbour1992}) for the approximation:} the binomial random variable $Z_n^{(l)}$ may be approximated by a Poisson random variable $\tilde{Z}_n^{(l)}\stackrel{(d)}{=}\mathcal{P}\big((\sum_{k=1}^{m}I_n^{(k)}\lambda_{kl})(\frac{N_l}{N}-\frac{A_{n-1}^{(l)}}{N}-\frac{U_n^{(l)}}{N})\big)$ such that
		\begin{align*}
		d_{\text{TV}}(Z^{(l)}_n, \tilde{Z}_n^{(l)}) \leq \frac{2}{(N_l - A_n^{(l)} - U_n^{(l)}) \left(\frac{\sum_{k=1}^{m}I_n^{(k)} \lambda_{kl}}{N}\right)} \sum_{i=1}^{N_l - A_n^{(l)} - U_n^{(l)}} \left(\frac{\sum_{k=1}^{m}I_n^{(k)}\lambda_{kl}}{N} \right)^2 \leq \frac{2\displaystyle\max_{k,l}\lambda_{kl}}{N}=\frac{2\lambda}{N}.
		\end{align*}
		As a consequence,  the first term in the right hand side of \eqref{EC_n} can be approximated  as
		\begin{equation}
		\left| \mathbb{E}\left[Z^{(l)}_n\mathbb{1}_{\sum_{j=1}^{m}Z^{(j)}_n\leq c}\bigg|\mathcal{F}_{n-1}\right]- \mathbb{E}\left[\tilde{Z}^{(l)}_n\mathbb{1}_{\sum_{j=1}^{m}\tilde{Z}^{(j)}_n\leq c} \bigg|\mathcal{F}_{n-1}\right] \right|\leq \frac{2mc\lambda}{N},
		\end{equation}
		and
		\begin{equation}
		\left|\mathbb{E}\left[\frac{Z_n^{(l)}}{\sum_{j=1}^{m}Z^{(j)}_n}\mathbb{1}_{\sum_{j=1}^{m} Z^{(j)}_n> c}\bigg|\mathcal{F}_{n-1}\right] - \mathbb{E}\left[\frac{\tilde{Z}_n^{(l)}}{\sum_{j=1}^{m}\tilde{Z}^{(j)}_n}\mathbb{1}_{\sum_{j=1}^{m}\tilde{Z}^{(j)}_n> c} \bigg|\mathcal{F}_{n-1}\right] \right|\leq \frac{2m \lambda}{N}.
		\end{equation}
		It follows that we need to deal with the Poisson random variables $\tilde{Z}^{(l)}_n (l \in \{1,...,m\})$. Because of the result that the sum of two independent Poisson random variables is a Poisson random variable whose parameter is the sum of the two parameters, we have that $\sum_{j\neq l}\tilde{Z}_n^{(j)}=: \hat{Z}_n^{(l)}$ {has a Poisson distribution} with parameter $\hat{\lambda}_n^{(l)}:= \sum_{j\neq l}\lambda_n^{(j)} $. And hence,
		\begin{align*}
		\mathbb{E}\bigg[\tilde{Z}^{(l)}_n\mathbb{1}_{\sum_{j=1}^{m}\tilde{Z}^{(j)}_n\leq c}\big|\mathcal{F}_{n-1}\bigg]&= \sum_{h=1}^{c} \sum_{h_1 = 1}^{h} h_1\frac{(\lambda_n^{(l)})^{h_1}(\hat{\lambda}_n^{(l)})^{h-h_1}}{h_1!(h-h_1)!}e^{-\Lambda_n} \\&= \lambda_n^{(l)}\sum_{h=1}^{c}  \frac{(\Lambda_n)^{h-1}}{(h-1)!}e^{-\Lambda_n} = \lambda_n^{(l)}\sum_{h=0}^{c} \frac{h}{\Lambda_n} \frac{(\Lambda_n)^{h}}{h!}e^{-\Lambda_n}
		\end{align*}
		and
		\begin{align}
		\mathbb{E}\bigg[\frac{\tilde{Z}_n^{(l)}}{\sum_{j=1}^{m}\tilde{Z}^{(j)}_n} \mathbb{1}_{\sum_{j=1}^{m}\tilde{Z}^{(j)}_n> c}\bigg|\mathcal{F}_{n-1}\bigg]&=\sum_{h=c+1}^{\infty}\sum_{k=0}^{h}\frac{k}{h}\frac{(\lambda^{(l)}_n)^k}{k!}\frac{(\hat{\lambda}_n^{(l)})^{h-k}}{(h-k)!}e^{-\lambda^{(l)}_n } e^{-\hat{\lambda}_n^{(l)}}\nonumber\\
		&=\lambda^{(l)}_n\sum_{h=c+1}^{\infty}\sum_{k=0}^{h-1}\frac{1}{h}\frac{(\lambda^{(l)}_n)^k}{k!}\frac{(\hat{\lambda}_n^{(l)})^{h-1-k}}{(h-1-k)!}e^{-\lambda^{(l)}_n } e^{-\hat{\lambda}_n^{(l)}}\nonumber\\
		&=\lambda^{(l)}_n\sum_{h=c+1}^{\infty}\frac{1}{h}\frac{(\lambda^{(l)}_n+ \hat{\lambda}_n^{(l)})^{h-1}}{(h-1)!} e^{-(\lambda^{(l)}_n+ \hat{\lambda}_n^{(l)})} =\frac{\lambda^{(l)}_n}{\Lambda_n}\sum_{h=c+1}^{\infty}\frac{(\Lambda_n)^{h}}{h!} e^{-\Lambda_n} \nonumber\\
		&= \frac{\lambda^{(l)}_n}{\Lambda_n}(1-\sum_{h=0}^{c}\frac{(\Lambda_n)^{h}}{h!}e^{-\Lambda_n} ).
		\end{align}
		{Using \eqref{EC_n}, we obtain:}
		\begin{equation*}
		\mathbb{E}[C_n^{(l)}|\mathcal{F}_{n-1}] = \mathbb{E}\left[\tilde{Z}^{(l)}_n\mathbb{1}_{\sum_{j=1}^{m}\tilde{Z}^{(j)}_n\leq c}  + \frac{\tilde{Z}_n^{(l)}}{\sum_{j=1}^{m}\tilde{Z}^{(j)}_n} \mathbb{1}_{\sum_{j=1}^{m}\tilde{Z}^{(j)}_n> c}\bigg|\mathcal{F}_{n-1}\right] = \frac{\lambda_n^{(l)}}{\Lambda_n}\left( c - \sum_{h=0}^{c} (c-h)  \frac{(\Lambda_n)^{h}}{h!}e^{-\Lambda_n} \right),
		\end{equation*}
		which finishes step 1.\\
	\textbf{Step 2:} We decompose the second term in the left hand side of \eqref{EC_n_in_lemma} as follow
		\begin{align} \label{decompositionofLambda}
		\frac{\lambda_n^{(l)}}{\Lambda_n}\left( c - \sum_{h=0}^{c} (c-h)  \frac{(\Lambda_n)^{h}}{h!}e^{-\Lambda_n} \right) =\alpha_n^{(l)}+\xi_n^{(l)}, \quad l = 1,...,m.
		\end{align}
		with
		\begin{align*}
		\alpha_n^{(l)}&:=\E\bigg[\frac{\lambda_n^{(l)}}{\Lambda_n}\left( c - \sum_{h=0}^{c} (c-h)  \frac{(\Lambda_n)^{h}}{h!}e^{-\Lambda_n} \right)\bigg|\mathcal{F}_{n-1}\bigg]\\
		\xi_n^{(l)}&:=  \frac{\lambda_n^{(l)}}{\Lambda_n}\left( c - \sum_{h=0}^{c} (c-h)  \frac{(\Lambda_n)^{h}}{h!}e^{-\Lambda_n} \right) -\E\bigg[\frac{\lambda_n^{(l)}}{\Lambda_n}\left( c - \sum_{h=0}^{c} (c-h)  \frac{(\Lambda_n)^{h}}{h!}e^{-\Lambda_n} \right)\bigg|\mathcal{F}_{n-1}\bigg].
		\end{align*}
		By writing
		\begin{align*}
		\alpha_n^{(l)} = \sum_{k=1}^{m}\P(I^{(k)}_n=1) \frac{\lambda_n^{k,l}}{\Lambda^k_n}\left( c - \sum_{h=0}^{c} (c-h)  \frac{(\Lambda^k_n)^{h}}{h!}e^{-\Lambda^k_n} \right),
		\end{align*}
		where
		\begin{equation}
		\lambda_n^{k,l}:= \lambda_{kl}\left(\frac{N_l}{N}-\frac{A_{n-1}^{(l)}}{N} - \frac{U_{n-1}^{(l)}}{N}- \frac{\mathbb{1}_{\{k=l\}}}{N}\right) \quad \text{and} \quad \Lambda^k_n:= \sum_{j=1}^{m}\lambda^{k,j}_n \quad (l = 1,...,m),
		\end{equation}
		we obtain that for every $t\in [0,1]$,
		\begin{equation}\label{reimannsumofalpha}
		\frac{1}{N}\sum_{n=1}^{\lfloor Nt \rfloor}\alpha_n^{(l)}= \frac{1}{N}\sum_{n=1}^{\lfloor Nt \rfloor} \left\{ \sum_{k=1}^{m} \dfrac{A_{n-1}^{(k)}}{|A_{n-1}|}\frac{\lambda_n^{k,l}}{\Lambda^k_n}\left( c - \sum_{h=0}^{c} (c-h)  \frac{(\Lambda^k_n)^{h}}{h!}e^{-\Lambda^k_n} \right)\right\}.
		\end{equation}
		It is obvious that $\frac{1}{N}\sum_{n=1}^{\lfloor Nt \rfloor}\xi_n $ is an $\mathcal{F}_t^N-$martigale with the quadratic variation,
		\begin{equation*}
		\langle \frac{1}{N}\sum_{n=1}^{\lfloor N\cdot \rfloor}\xi_n \rangle_t = \frac{1}{N^2}\sum_{n=1}^{\lfloor Nt \rfloor}\E\left[\xi_n^2|\mathcal{F}_{n-1}\right] \leq \frac{1}{N^2}\sum_{n=1}^{\lfloor Nt \rfloor} m(c+1)^2 \leq \frac{m(c+1)^2}{N}.
		\end{equation*}
		By the Doob's inequality, we have
		\begin{equation*}
		\E\left[ \max_{0\leq s \leq t}\|\frac{1}{N}\sum_{n=1}^{\lfloor Nt \rfloor}\xi_n \|^2\right] \leq 4 \E\left[\| \langle \frac{1}{N}\sum_{n=1}^{\lfloor N\cdot \rfloor}\xi_n \rangle_t \|\right]\leq \frac{4m(c+1)^2}{N} \stackrel{N\rightarrow \infty}{\longrightarrow} 0,
		\end{equation*}
		which deduces that as N tends to infinity, we have that
		\begin{equation}\label{limitofxi_n}
		\frac{1}{N}\sum_{n=1}^{\lfloor Nt \rfloor}\xi_n \stackrel{(\mathbb{L}^2)}{\rightarrow} 0
		\end{equation}
		uniformly in $t\in[0,1]$.
		{Together with the points given in \eqref{EC_n_in_lemma}, \eqref{decompositionofLambda} and \eqref{limitofxi_n}, take the limit as $N \rightarrow \infty$ in the right hand side of \eqref{reimannsumofalpha}, we obtain the right hand side of \eqref{limitofC_n}.}\\
		\textbf{Step 3:} We use similar arguments as in step 2 to obtain the limit in right hand side of \eqref{lim:limitofHn}.
		Denote by $$\mu_n^{(l)}:=\left(\sum_{k=1}^{m}I_n^{(k)}\lambda_{kl}\right)\left(\frac{N_l}{N}-\frac{A_{n-1}^{(l)}}{N}-\frac{B_{n-1}^{(l)}}{N}-\frac{U_n^{(l)}}{N}\right).$$
		Recall from \eqref{distributionofHn} that conditioning on $\mathcal{F}_{n-1}$, $H_n^{(l)} \stackrel{(d)}{=} \mathcal{B}\textit{in} \left(N_l-A_{n-1}^{(l)}-B_{n-1}^{(l)}-U_n^{(l)}, \frac{\sum_{k=1}^{m}I_n^{(k)} \lambda_{kl}}{N}\right)$, then
		\begin{equation}\label{EH_n}
		\frac{1}{N}\sum_{n=1}^{\lfloor Nt \rfloor}\E[H_n^{(l)}|\mathcal{F}_{n-1}] = \frac{1}{N}\sum_{n=1}^{\lfloor Nt \rfloor} \mu_n^{(l)}.
		\end{equation}
	We write
		\begin{equation}
		\frac{1}{N}\sum_{n=1}^{\lfloor Nt \rfloor}\mu_n^{(l)} = \frac{1}{N}\sum_{n=1}^{\lfloor Nt \rfloor}(\beta_n^{(l)} + \zeta_n^{(l)})
		\end{equation}
		with
		\begin{align*}
		\beta_n^{(l)}&:= \E\left[\left(\sum_{k=1}^{m}I_n^{(k)}\lambda_{kl}\right)\left( \frac{N_l}{N}-\frac{A_{n-1}^{(l)}}{N}-\frac{B_{n-1}^{(l)}}{N}-\frac{U_n^{(l)}}{N}\right) \bigg| \mathcal{F}_{n-1} \right];\\
		\zeta_n^{(l)}&:= \left(\sum_{k=1}^{m}I_n^{(k)}\lambda_{kl}\right)\left(\frac{N_l}{N}-\frac{A_{n-1}^{(l)}}{N}-\frac{B_{n-1}^{(l)}}{N}-\frac{U_n^{(l)}}{N}\right) - \E\left[\left(\sum_{k=1}^{m}I_n^{(k)}\lambda_{kl}\right)\left(\frac{N_l}{N}-\frac{A_{n-1}^{(l)}}{N}-\frac{B_{n-1}^{(l)}}{N}-\frac{U_n^{(l)}}{N}\right) \bigg| \mathcal{F}_{n-1}\right].
		\end{align*}
		Using a similar argument as in step 2, we have
		\begin{align*}
		\frac{1}{N}\sum_{n=1}^{\lfloor Nt \rfloor}\beta_n^{(l)}& =\frac{1}{N}\sum_{n=1}^{\lfloor Nt \rfloor}\sum_{k=1}^{m} \mathbb{P}(I_n^{(k)} =1) \lambda_{kl}\left( \frac{N_l}{N}- \frac{A^{(l)}_{n-1}}{N} - \frac{B^{(l)}_{n-1}}{N} - \frac{U^{(l)}_{n-1}}{N} -\frac{\mathbb{1}_{\{k\neq l\}}}{N}\right)\\
		&= \frac{1}{N}\sum_{n=1}^{\lfloor Nt \rfloor}\sum_{k=1}^{m}\dfrac{A_{n-1}^{(k)}}{\|A_{n-1}\|}  \mu_n^{k,l} - \frac{1}{N}\sum_{n=1}^{\lfloor Nt \rfloor}\sum_{k=1}^{m} \dfrac{A_{n-1}^{N,k}}{\|A_{n-1}^N\|} \lambda_{kl} \frac{\mathbb{1}_{\{k\neq l\}}}{N},
		\end{align*}
		with $\mu^{k,l}_n := \lambda_{kl}\left( \frac{N_l}{N}- \frac{A^{(l)}_{n-1}}{N} - \frac{B^{(l)}_{n-1}}{N} - \frac{U^{(l)}_{n-1}}{N}\right)$. Then,
		\begin{equation}
		    \left| \frac{1}{N}\sum_{n=1}^{\lfloor Nt \rfloor} \left(\beta_n^{(l)} -  \sum_{k=1}^{m}\dfrac{A_{n-1}^{(k)}}{\|A_{n-1}\|}  \mu_n^{k,l} \right) \right| \leq \frac{1}{N}\sum_{n=1}^{\lfloor Nt \rfloor}\sum_{k=1}^{m} \dfrac{A_{n-1}^{N,k}}{\|A_{n-1}^N\|} \lambda_{kl} \frac{\mathbb{1}_{\{k\neq l\}}}{N} \leq \frac{\lambda}{N}.
		\end{equation}
	{Take the limit as $N \rightarrow + \infty$, we have that \[\lim_{N \rightarrow +\infty} \frac{1}{N}\sum_{n=1}^{\lfloor Nt \rfloor} \sum_{k=1}^{m}\dfrac{A_{n-1}^{(k)}}{\|A_{n-1}\|}  \mu_n^{k,l} = \int_{0}^{t} \sum_{k=1}^m \frac{a_s^{(k)}}{\|a_s\|} \mu_s^{k,l}ds.\]
	Further, the $\mathcal{F}_{t}^N-$martingale $\displaystyle\frac{1}{N}\sum_{n=1}^{\lfloor N\cdot \rfloor}\zeta_n^{(l)}$ {converges in $\mathbb{L}^2$ to $0$ uniformly in $t \in [0,1]$}. Thus, \eqref{lim:limitofHn} is proved.}\\
	{For the proof of \eqref{lim:limitofIn}, by the definition of $I_n$ as in \eqref{eq:defofIn}, we have
	\[ \frac{1}{N}\sum_{n=1}^{\lfloor Nt \rfloor} \E[I_n^{(l)}|\mathcal{F}_{n-1}] =\frac{1}{N}\sum_{n=1}^{\lfloor Nt \rfloor}  \frac{A_{n-1}^{(l)}}{\|A_{n-1}\|} = \frac{1}{N}\sum_{n=1}^{\lfloor Nt \rfloor}  \frac{A_{n-1}^{(l)}/N}{\|A_{n-1}\|/N} .
	\] Take the limit as $N \rightarrow +\infty$, we obtain the limit in the right hand side of \eqref{lim:limitofIn}.}\\
	
	The preceding steps allow to conclude the proof of Proposition \ref{prop:conv-Delta}.
\end{proof}	
	
	\subsection{The uniqueness}\label{sec:uniqueness}
	It remains to prove that the {limiting value $x=(a,b,u)$ we have found is unique solution of the system of he ODEs \eqref{limiting_eq}. If it is the case, then the process $(X^N)_N$ admits a unique limiting value and thus converges to $x$}.\\
	Assume that there exist two solutions $x^1$ and $x^2$  to ODEs \eqref{limiting_eq} on the interval $[0,t'_0]$, where
	\[ t'_0 = \inf\{ t \in [0,1]: a^1_{t'_0} =0 \text{ or } a^2_{t'_0} = 0 \}.\] Then using the intermediate value theorem, there exist $ \xi_{ij}(s) \in [x^1_{ij}(s), x^2_{ij}(s)]$ such that
	\begin{align*}
	\|x^{1}_t - x^2_t\| = \left\|\int_{0}^{t}(f(x^1_s) - f(x_s^2))ds\right\| &\leq \int_{0}^{t}\sum_{i=1}^{3}\sum_{j=1}^{m} \left|\frac{\partial f}{\partial x_{ij}}(\xi_{ij}(s)) \right| \big|x_{ij}^1(s)-x_{ij}^2(s)\big|ds\\
	&\leq \int_{0}^{t}L(s)\|x^1_s-x^2_s\|ds,
	\end{align*}
	where $x^k_s=(x_{ij}(s))_{\substack{1\leq i \leq 3 \\1\leq j \leq m}}$ ($k\in \{1,2\}$) and $L(s)= \displaystyle\sum_{i=1}^{3}\sum_{j=1}^{m} \max \left|\frac{\partial f}{\partial x_{ij}}(x_s) \right|$, of which the maximum is over $x(s)= (x_{ij}(s))_{ij} \in [0,1]^{3m} $ such that $\forall i,j: x_{ij} \in [x_{ij}^1, x^2_{ij}]$, where by an abuse of notation, the bounds of interval $[x^1_{ij}, x^2_{ij}]$ can be switched depending on the minimum or maximum of the bounds.\\
	 By the Gr\"onwall's inequality, we get
	\begin{align*}
	\|x^{1}_t - x^2_t\| \leq \|x^{1}_0 - x^2_0\| \exp(\int_{0}^{t}L(s)ds)=0.
	\end{align*}
	This shows that $x^1_t \equiv x^2_t$ for all $t \in [0,t'_0]$. It also follows $t'_0 =t_0$.
	
	\section{Simulation}
	
\begin{figure}
  	\begin{tabular}{c c}
  \begin{subfigure}{.5\textwidth}
  	\centering
  	\includegraphics[width=.85\linewidth]{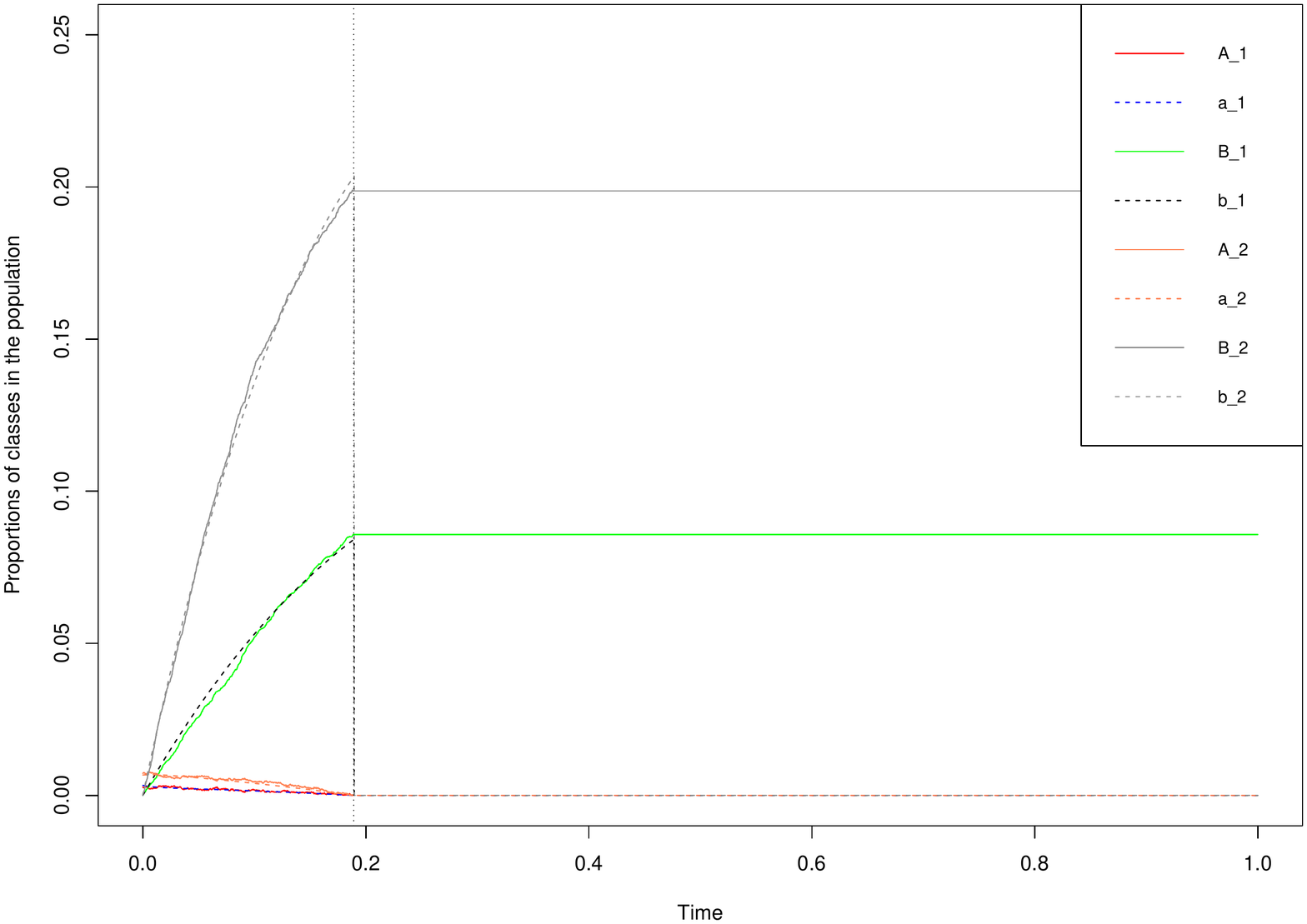}
  	\caption{$c=1$}
  	\label{fig:sfig1}
  \end{subfigure} &
\begin{subfigure}{.5\textwidth}
	\centering
	\includegraphics[width=.85\linewidth]{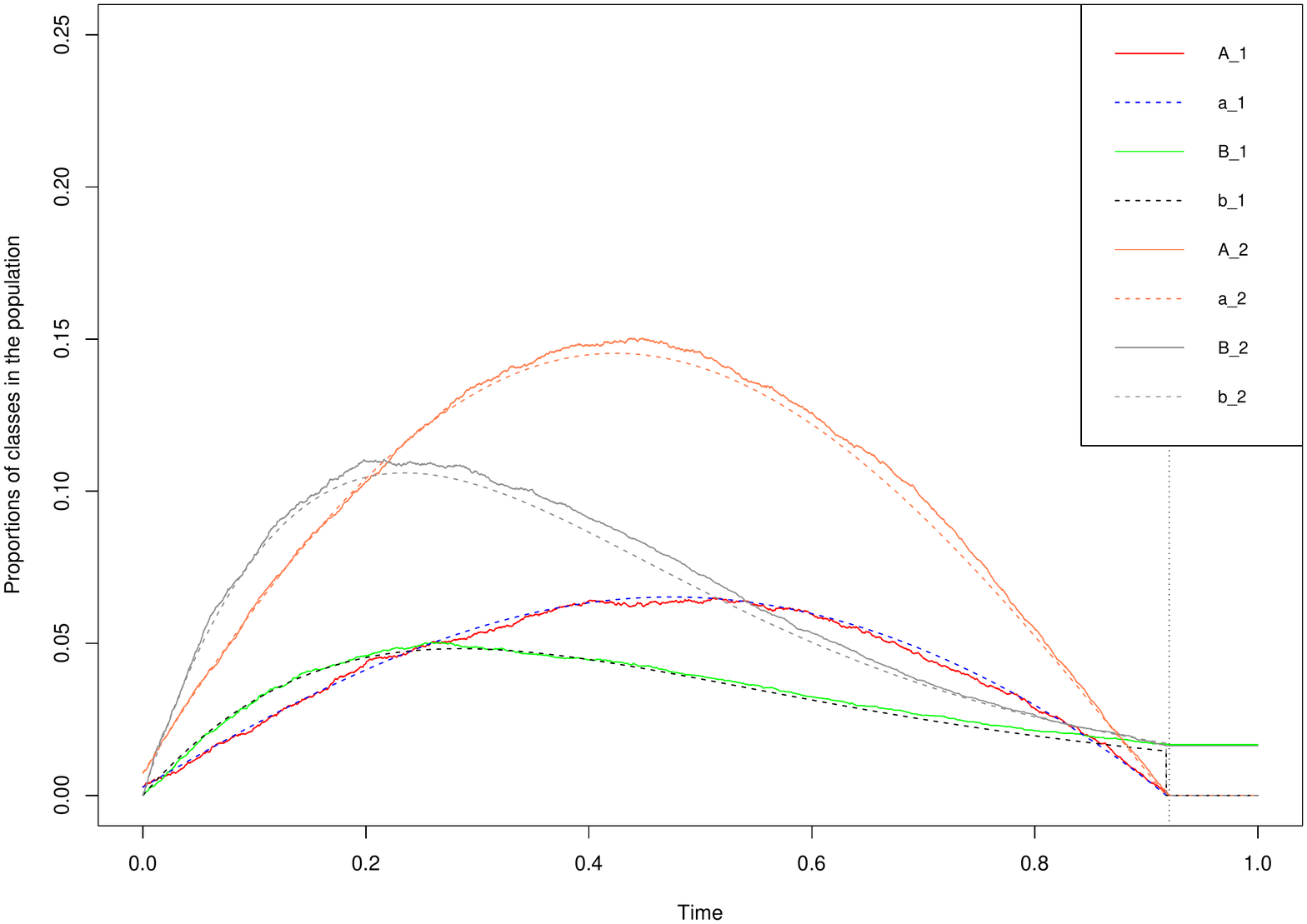}
	\caption{$c=2$}
	\label{fig:sfig2}
\end{subfigure}\\

\begin{subfigure}{.5\textwidth}
	\centering
	\includegraphics[width=.85\linewidth]{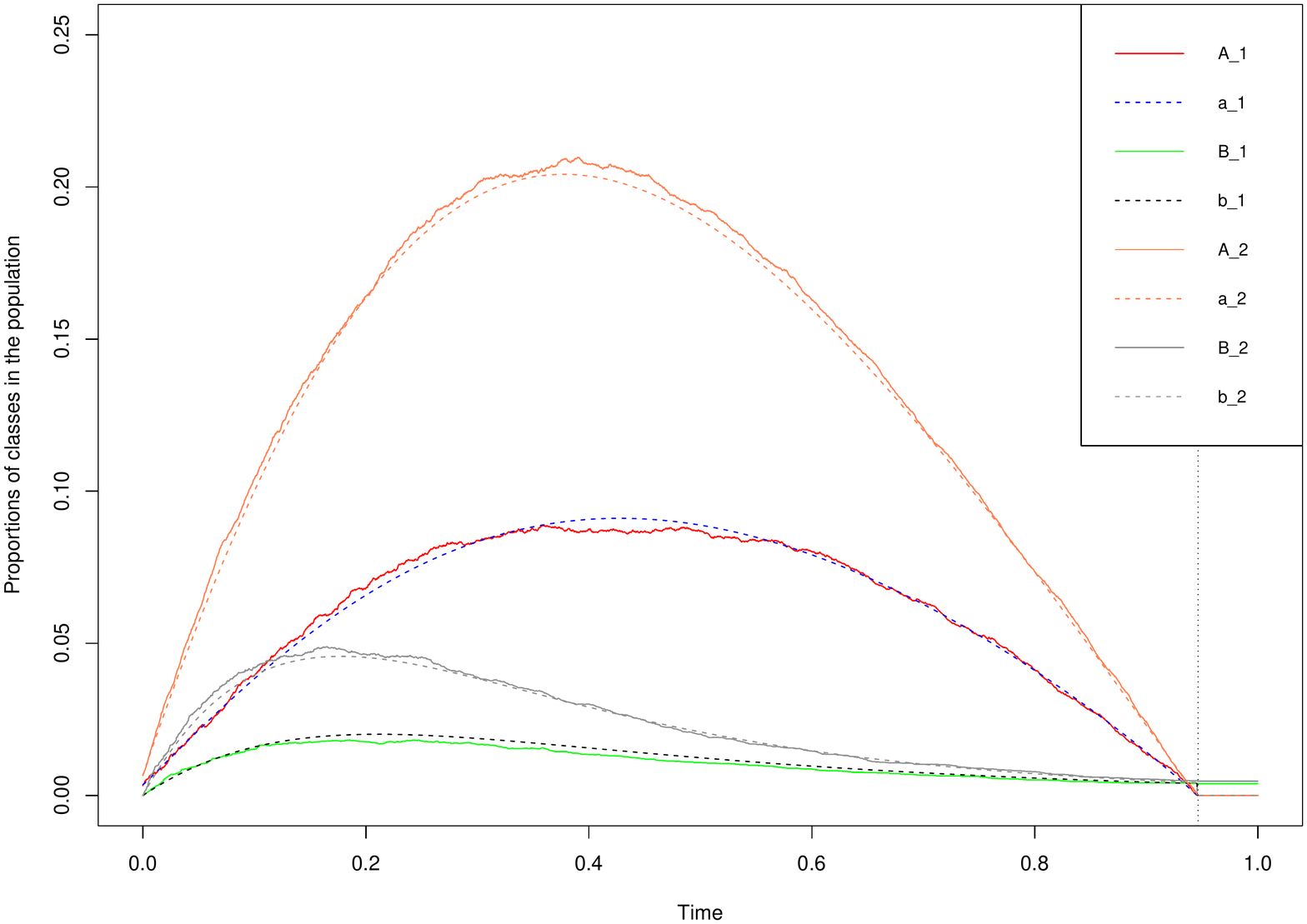}
	\caption{$c=3$}
	\label{fig:sfig3}
\end{subfigure} &
\begin{subfigure}{.5\textwidth}
	\centering
	\includegraphics[width=.85\linewidth]{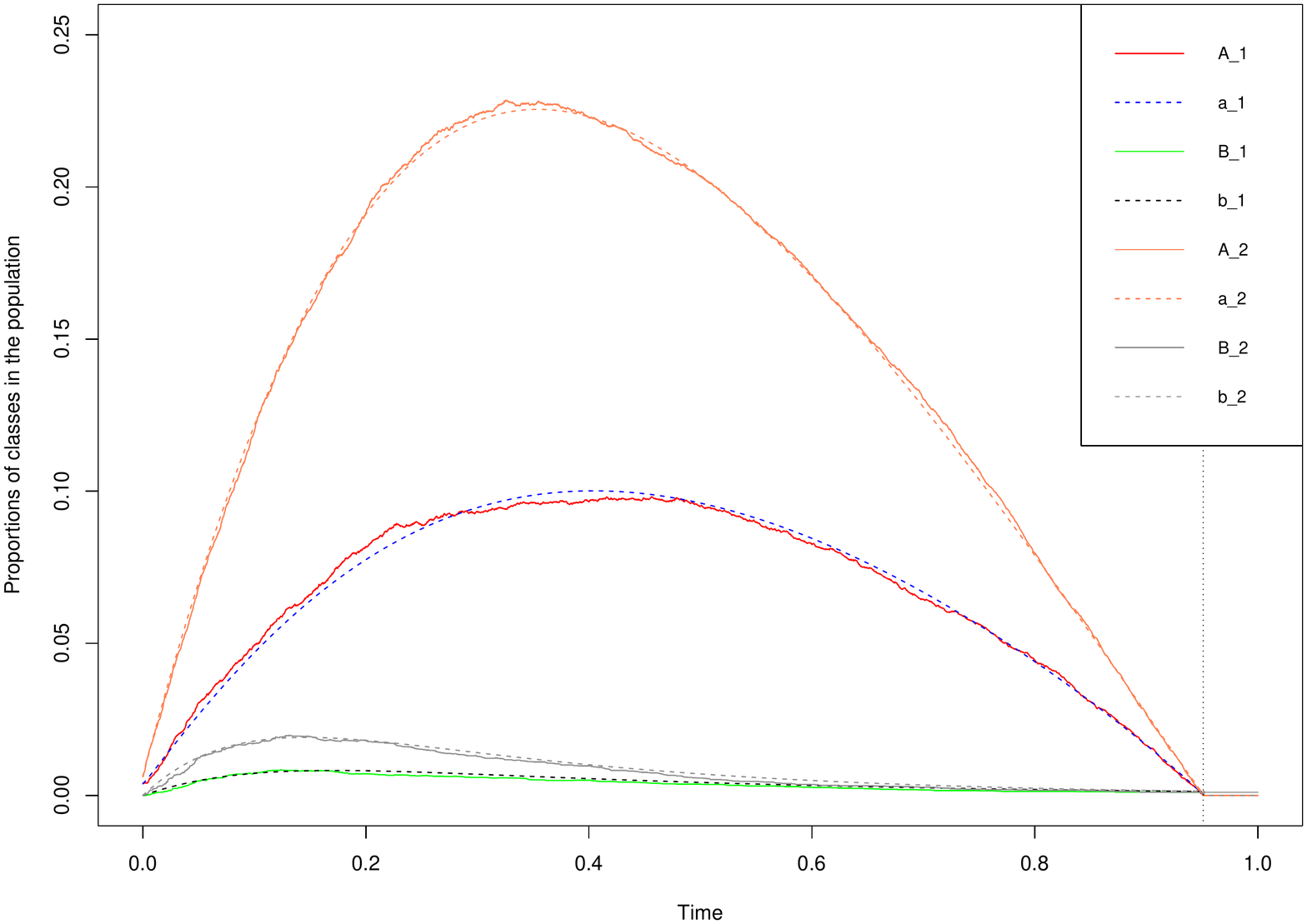}
	\caption{$c=4$}
	\label{fig:sfig4}
\end{subfigure}\\
\begin{subfigure}{.5\textwidth}
	\centering
	\includegraphics[width=.85\linewidth]{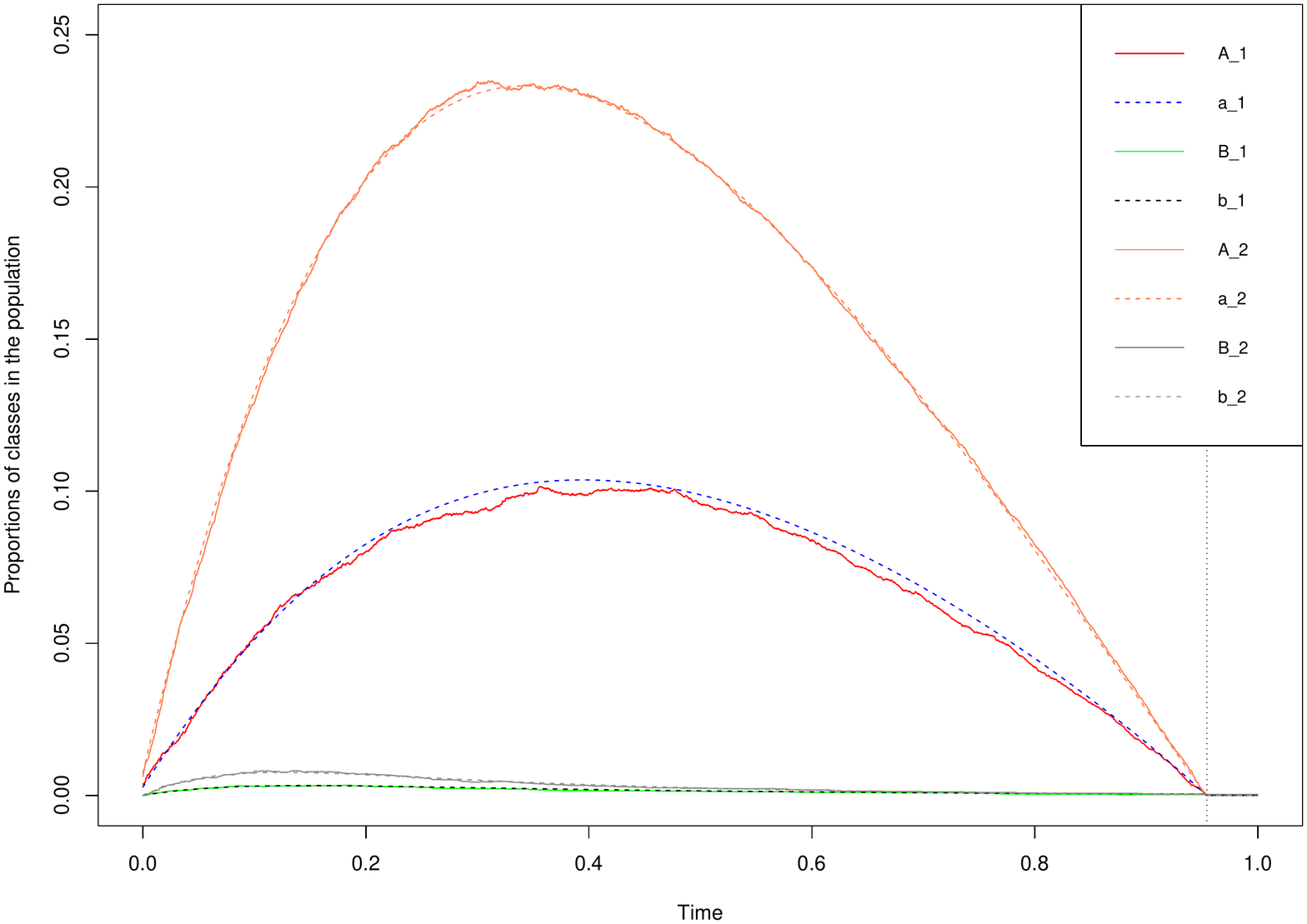}
	\caption{$c=5$}
	\label{fig:sfig5}
\end{subfigure}&
\begin{subfigure}{.5\textwidth}
	\centering
	\includegraphics[width=.85\linewidth]{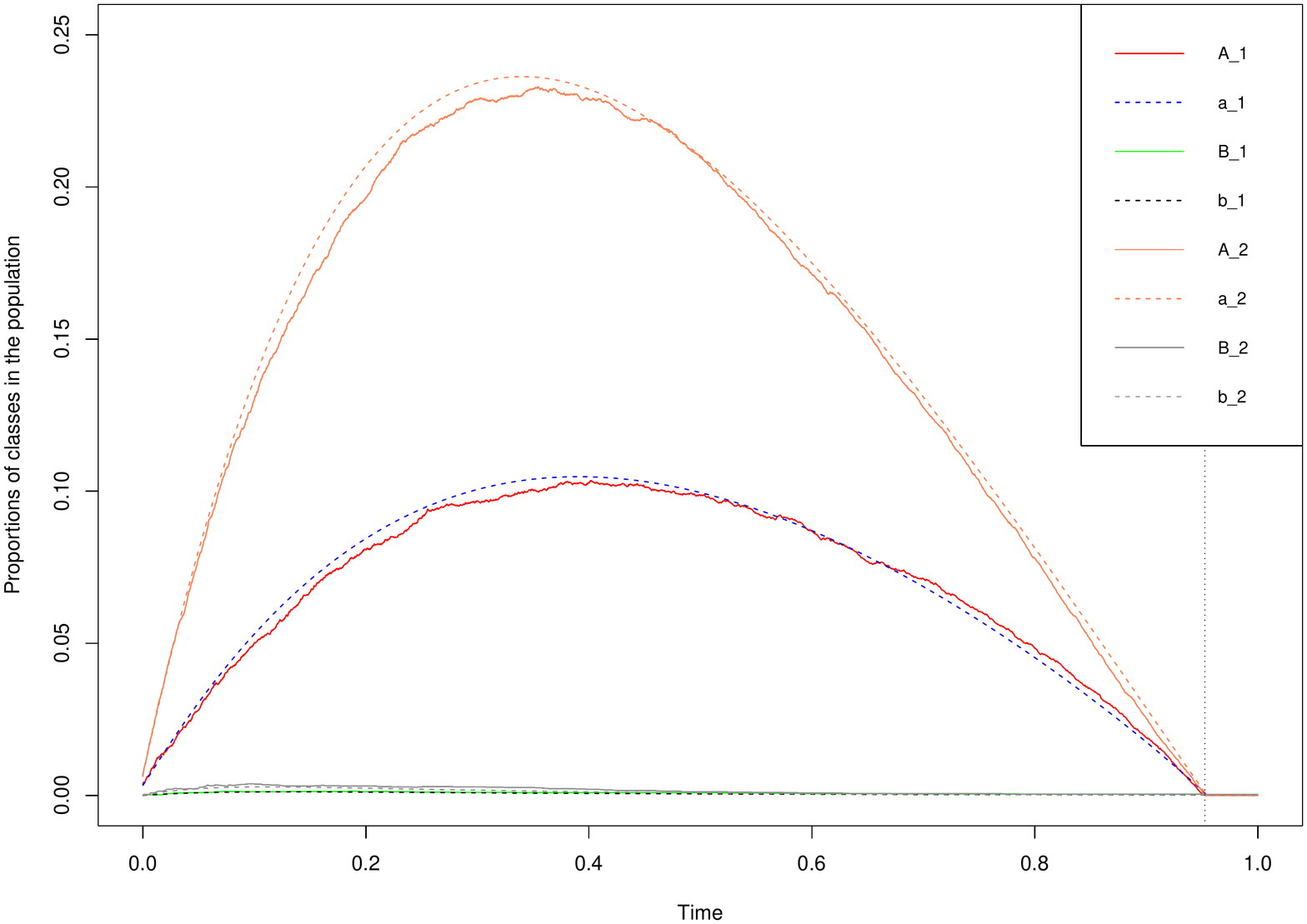}
	\caption{$c=6$}
	\label{fig:sfig5}
\end{subfigure}
  	\end{tabular}

\caption{Plots of the proportions of classes in the population of size $N=10000$ when  $c$ varies from $1$ to $6 $ and all the others parameters are fixed: $\|A_0\|=100$ the parameters $\pi=(1/3,2/3)$, $\lambda_{11}=2,\lambda_{12}=3, \lambda_{22}=4$.}
\label{fig:AnBn-Odes}
\end{figure}

\begin{figure}[hbtp!]
	\centering
	\includegraphics[width=8cm,height=8cm]{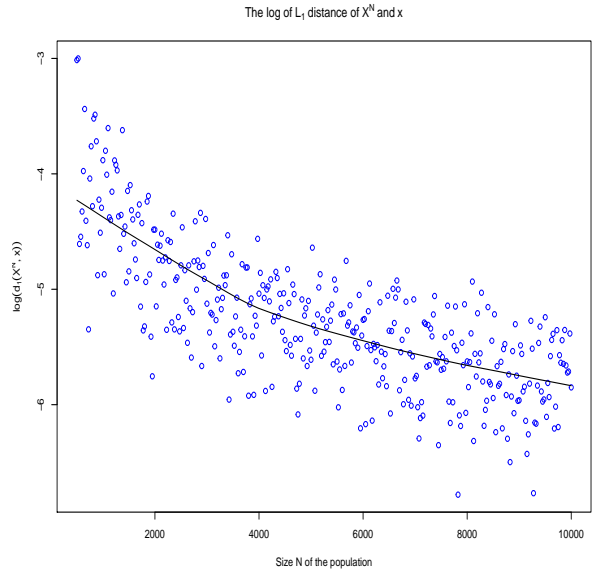}
	\caption{Scatter plot of $\ln d_1(X^N,x)$ along with the smoothing line suggesting the linear relationship between $\ln d_1(X^N,x)$ and $N$. The plot is done for the case $c=3$, the number of initial individuals are $1\%$ of the population and the size $N$ varies from $500$ to $10000$. All other parameters are fixed: $\pi=(1/3,2/3)$, $\lambda_{11}=2,\lambda_{12}=3, \lambda_{22}=4$.}
	\label{fig:convergence}
\end{figure}

	The simulations show that the deterministic solution of the system of ODEs \eqref{limiting_eq} fits well with our stochastic process, see Figure \ref{fig:AnBn-Odes}. {The sequence of stochastic process $(X_\cdot^N)_N$ that we have constructed describes how the chain-referral process works on a network. When we consider the population with a very large number of people, the process $(X_\cdot^N)_N$ is asymptotically  a deterministic function, which is a solution of a system of \eqref{limiting_eq}. To see numerically the convergence, we do a simulation: for $c=3$, we vary $N$ from 500 to 50000 and plot as a function of $N$ the $\log$ of the quantity:
	\[ \int_0^1 (\|A^N_t-a_t\| + \|B^N_t-b_t\| +\|U^N_t - u_t\| ) dt,\]
	Figure \ref{fig:convergence}. The speed of convergence has been studied in the case of Erdös-Rényi graphs in the PhD-thesis, by establishing a central limit theorem \cite{vothesis}.}\\

	 By studying the solution of \eqref{limiting_eq}, we can obtain an approximation of the fraction of the population that has been interviewed {when the CRS process stops}. {The proportion of the population discovered is then approximated by $t_0$.}
	
	The number of maximum coupon $c$ plays an important role in how many people we could explore before there is no distributed coupons any more (when $\|a_t\|=0$). By keeping all other parameters fixed and changing $c$, in the simulations of Figure \ref{fig:AnBn-Odes}, we see that the time $t_0$ are different. For example, with $m=2$, $\pi=(1/3,2/3)$, $\lambda_{11}=2, \lambda_{22}=4, \lambda_{12}=3$, we obtain the table \ref{tab1}.
	\begin{table}[h!]
		\centering
		\begin{tabular}{|c|c|c|c|c|c|c|c|}
			\hline
			c& 1& 2& 3 & 4  &5 & 6 & ...\\
			\hline
			$t_0$& 0.18&0.91&0.94& 0.95& 0.95 & 0.95 &...\\
			\hline
		\end{tabular}
	\caption{Numerical computation of $t_0$ for varying parameters $c$.} \label{tab1}
	\end{table}

\begin{figure}[hbtp!]
		\centering
		\includegraphics[width=10cm,height=10cm]{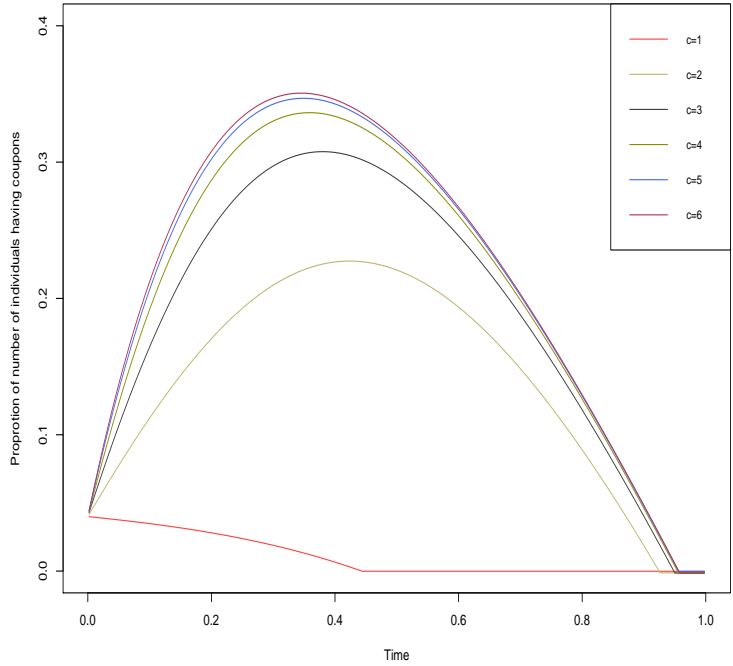}
		\caption{Plot the function $\|a\|$ for $6$ cases: $c$ takes values from $1$ to $6$. All other parameters are fixed: $\|a_0\|=0.05, \pi=(1/3,2/3)$, $\lambda_{11}=2,\lambda_{12}=3, \lambda_{22}=4$. The values $\|a_t\|$ represents the proportion of individuals having coupons at time $t$.}
		\label{fig:totalcoupons}
	\end{figure}

	If $c=1$, even though the average number of neighbors are bigger than $1$, the simple random walk describing the survey reaches only a very small number of people, see Figure \ref{fig:sfig1}. The random walk stops when it encounters a node of degree 1 and can not propagate any more. \\
	Furthermore, the parameter $c$ also impacts the peaks (time and size) the curves corresponding to the number of distributed coupons. In case of a limited budget with a fixed number of interviews, a higher value of $c$ can imply that we discover a larger fraction of the population since it allows more flexibility in the interviewees. From the Figure \ref{fig:totalcoupons}, we observe that the proportion of people receiving coupons gets bigger as $c$ increases. If $c=1$, the fraction of discovered population is small, which means that the survey is not so efficient. When $c$ takes values from $4$ to $6$, the corresponding curves of $\|a_t\|$ are "close" and so are the times $t_0$. However, in these cases, the number of coupons spent during the CRS survey is large. 
		{We can also be interested in seeing how $c$ impacts the part of population discovered when the survey stops after a fixed number of interviewed individuals. For example, consider the case when $N= 1000$ and assume that we start with $A_0 = 10$. The parameters of the SBM are $\pi=(1/3,2/3)$, $\lambda_{11}=2, \lambda_{22}=4, \lambda_{12}=3$. Then when there have been approximately $\lfloor 0.2N \rfloor$ individuals interviewed, the proportion of the explored individuals: $\|A^N_{0.2} \| + \|B^N_{0.2}\|$ for each $c$ varying from $1$ to $6$ is given in Table \ref{tab2}}.\\
\begin{table}[h!]
	\centering
	\begin{tabular}{|c|c|c|c|c|c|c|c|}
		\hline
		c& 1& 2& 3 & 4  &5 & 6 \\
		\hline
		$\|A^{1000}_{0.2}\| + \|B^{1000}_{0.2}\|$& 0.213 & 0.308 &0.268 &0.308 &0.310 &0.260\\
		\hline
	\end{tabular}
	\caption{Numerical computation of $\|A^N_{t}\| + \|B^N_{t}\|$ for varying parameters $c \in \{1, \ldots, 6\}$ at time $t=0.2$ and $N=1000, A_0 =10,\pi=(1/3,2/3)$, $\lambda_{11}=2, \lambda_{22}=4, \lambda_{12}=3.$} \label{tab2}
\end{table}

\begin{figure}[hbtp!]
		\centering
		\includegraphics[width=8cm,height=8cm]{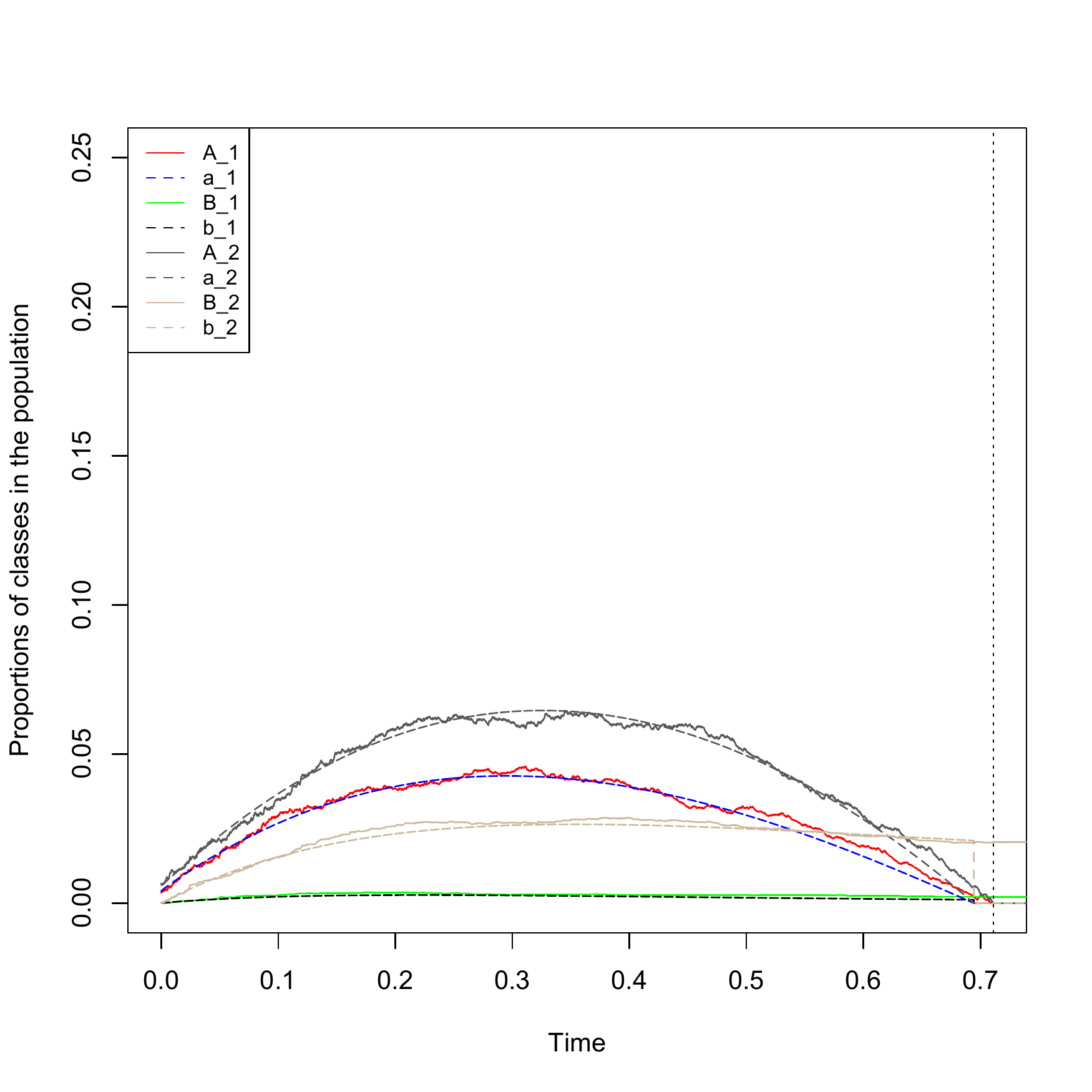}
		\caption{Plot the proportion of classes in the case $c=3, N=1000, A_0=10, \pi=(1/3,2/3)$ and the graph is bipartite $\lambda_{11}=\lambda_{22}=0, \lambda_{12} = 4$.}
		\label{fig:AnBn-lambda=0-4}
	\end{figure}

	Changing the parameters $\lambda_{kl}$ impacts the discovered proportion of types. For instant, let us take a bipartite random model $\pi = (1/3,2/3), c=3$ and $\lambda_{11}= \lambda_{22} = 0$, $\lambda_{12}=4$, which means that the people between communities are highly connected and there is no connection within community. In this case, the number of explored people without coupon of type 1 is quite small compared to the one of type 2, see Figure \ref{fig:AnBn-lambda=0-4}. 
	
	\clearpage


		
		\providecommand{\noopsort}[1]{}\providecommand{\noopsort}[1]{}\providecommand{\noopsort}[1]{}\providecommand{\noopsort}[1]{}\providecommand{\noopsort}[1]{}

\end{document}